\RequirePackage{fix-cm}
\documentclass[smallextended]{svjour3}       
\smartqed  
\usepackage{graphicx}
\usepackage{amssymb}
\usepackage[colorlinks=true,linkcolor=blue]{hyperref}
\usepackage{amsmath,amssymb}
\usepackage{graphicx,graphics}
\usepackage{subfigure}
\usepackage{dcolumn}
\usepackage{bm}
\usepackage{float}

\newcommand{\cQ}{\mathcal{Q}}

\renewcommand{\div}{\mbox{\rm div\,}}
\newcommand{\curl}{\mbox{\rm curl\,}}

\newcommand{\cX}{\mathcal{X}}

\newcommand{\cM}{\mathcal{M}}

\newcommand{\cW}{\mathcal{W}}

\newcommand{\R}{\mathbf{R}}

\begin{document}

\title{Error analysis of a fully discrete scheme for the
Cahn-Hilliard-Magneto-hydrodynamics problem
}

\titlerunning{Error analysis for the
Cahn-Hilliard-Magneto-hydrodynamics problem}        

\author{Hailong Qiu
}


\institute{Hailong Qiu \at
School of Mathematics and Physics, Yancheng Institute of Technology,
Yancheng, 224051, China.\\
 This work is supported by the Natural Science Foundation of China (No. 11701498).\\
              \email{qhllf@163.com}           
}

\date{Received: date / Accepted: date}

\maketitle

\begin{abstract}
In this paper we analyze a fully discrete scheme for
a general Cahn-Hilliard equation coupled with a nonsteady Magneto-hydrodynamics flow, which describes two immiscible, incompressible
and electrically conducting fluids with different mobilities, fluid viscosities and magnetic diffusivities.
A typical fully discrete scheme, which is comprised of
conforming finite element method and the Euler semi-implicit discretization
based on a convex splitting of the energy of the equation is considered in detail.
 We prove that our scheme is
unconditionally energy stability and obtain some optimal
 error estimates for the concentration field, the chemical potential, the velocity field, the magnetic field
and the pressure. The results of numerical tests are presented to
validate the rates of convergence.
\end{abstract}

\section{Introduction}
\label{intro}

In the paper, we derive error estimates for a fully discrete, first order in time, finite element method
for the Cahn-Hilliard-Magneto-hydrodynamics problem for two phase flow. Let $\Omega \subset \R^d$ ($d=2,3$)  stands for an
open convex polygonal or polyhedral domain with Lipschitz continuous boundary $\partial\Omega$. For all $\phi \in H^1(\Omega)$, $\textbf{u} \in \textbf{L}^2(\Omega)$
and $\textbf{B} \in \textbf{L}^2(\Omega)$, consider the following energy:
\begin{align}\label{eq1}
\textbf{E}(\phi,\textbf{u},\textbf{B})=\int_{\Omega}\Bigl(\frac{1}{2}\|\textbf{u}\|^2&+\frac{S_c}{2}\|\textbf{B}\|^2
+\frac{\lambda\epsilon}{2}\|\nabla\phi\|^2+\frac{\lambda}{4\epsilon}(1-\phi^2)^2\Bigr)dx,
\end{align}
where $\phi$, $\textbf{u}$ and $\textbf{B}$ denote respectively the concentration field,
the velocity field and the magnetic field, and the parameter $\epsilon>0$ stands for the interfacial thickness between the two phases.
The Cahn-Hilliard-Magneto-hydrodynamics system is a gradient flows of this energy:
\begin{subequations}\label{eq2.1}
\begin{alignat}{2} \label{eq2.1a} \partial_t \phi-\epsilon\mbox{div}(\kappa(\phi)\nabla\mu)-\nabla\phi\cdot\textbf{u}&=0,  &&\ \mbox{in}\  \Omega\times T,\\
-\epsilon\Delta \phi+\epsilon^{-1}(\phi^3-\phi)&=\mu,  &&\ \mbox{in}\  \Omega\times T,\\
\partial_t \textbf{u}-\mbox{div}(2\nu(\phi)\mathbb{D}(\textbf{u}))
+(\textbf{u}\cdot \nabla)\textbf{u}
+S_c\textbf{B}\times\mbox{curl}\textbf{B}
\\\qquad\qquad\qquad\qquad\nonumber+\nabla p&=\lambda\mu\nabla\phi,  &&\ \mbox{in}\  \Omega\times T,\\
\partial_t \textbf{B}+\mbox{curl}(\eta(\phi)\mbox{curl} \textbf{B})-\mbox{curl}(\textbf{u}\times \textbf{B})&=\textbf{0}, &&\ \mbox{in}\  \Omega \times T,\\
\mbox{div} \textbf{u}&=0, &&\ \mbox{in}\  \Omega\times T,\\
\mbox{div} \textbf{B}&=0, &&\ \mbox{in} \ \Omega\times T.
\end{alignat}
\end{subequations}
This problem is considered in conjunction
with natural and no-flux/no-flow boundary conditions:
\begin{subequations}\label{eq3.1}
\begin{alignat}{2} \label{eq3.1a}
\frac{\partial\phi}{\partial\textbf{n}}=\frac{\partial\phi}{\partial\textbf{n}}&=0,\quad &&\ \mbox{on} \ \partial\Omega \times T,\\
 \textbf{u}&=\textbf{0}, \quad &&\ \mbox{on} \ \partial\Omega \times T,\\
 \textbf{B}\cdot \textbf{n}&=0, \quad &&\ \mbox{on} \ \partial\Omega \times T,\\
\mbox{curl \textbf{B}}\times \textbf{n}&=\textbf{0}, \quad &&\ \mbox{on} \  \partial\Omega\times T,
\end{alignat}
\end{subequations}
and initial conditions:
\begin{align}\label{eq4}
\phi(\textbf{x},0)=\phi_0(\textbf{x}),\ \
\textbf{u}(\textbf{x},0)=\textbf{u}_0(\textbf{x}),\ \
\textbf{B}(\textbf{x},0)=\textbf{B}_0(\textbf{x}), \  \mbox{in}\ \ \Omega,
\end{align}
 where $T>0$ is time, $\phi\approx\pm1$ represents two different fluids,
 $\mu$ and $p$  are respectively the chemical potential
 and the pressure.
 $\mathbb{D}(\textbf{u}):=\frac{1}{2}(\nabla \textbf{u}+\nabla \textbf{u}^T)$ is the strain-rate tensor.
  The given functions $\kappa(\phi)>0$, $\nu(\phi)>0$ and $\eta(\phi)=\frac{1}{\eta_0\sigma(\phi)}>0$
stand for respectively the mobility, the fluid viscous,
 the magnetic diffusivity with $\sigma$ the
electrical conductivity. $\eta_0$  is the magnetic permeability.
$S_c:= \frac{1}{\eta_0\rho_0}$ denotes the coupling coefficient,  here $\rho_0$ is the reference density.
$\lambda$ denotes the mixing energy density.

The \eqref{eq2.1}-\eqref{eq4} system satisfies the following energy laws:
\begin{align}\label{eq5}
\textbf{E}(\phi(t),\textbf{u}(t),\textbf{B}(t))&+\int_{\Omega}\Bigl(\nu(\phi)\|\nabla\textbf{u}\|^2+S_c\eta(\phi)\|\nabla\textbf{B}\|^2
+\lambda\epsilon\kappa(\phi)\|\nabla\mu\|^2\Bigr)dx\\&\nonumber\quad=\textbf{E}(\phi_0,\textbf{u}_0,\textbf{B}_0).
\end{align}

In the last decades, phase field approaches for two phases incompressible flows have been widely developed
to model and numerically solve the topological transitions of interfaces (cf. \cite{Elliott1989,1996Gurtin,Anderson1998,Lowengrub1998,Jacqmin1999,L2003,Yue2004,2004Feng} and references therein).
Recently, the research of interaction of electromagnetic fields with two immiscible, incompressible and electrically conducting
 fluids has become more and more important
for the design and analysis of engineering field, such as fusion reactors,
pump accelerators, metallurgical industry and Magneto-hydrodynamics generators \cite{Szekely1979,Moreau1990,Morley2000,2001Davidson}. 
Very recently, in the work \cite{Y2019}, two phase Magneto-hydrodynamics problem
about the diffuse interface between two different incompressible fluids is considered.

Though many error estimates are available for fully discrete scheme to
the Magneto-hydrodynamics equations \cite{Wiedmer2000,Prohl2008,Tone2009,L2013,He2015}
and for the Cahn-Hilliard/Navier-Stokes equations \cite{Elliott1989,2004Feng,2006Feng,2007Feng,Diegel2015,Shen2015,2016Feng,Caishen2018},
 it is not trivial to cope with the system which couple Magneto-hydrodynamics
 with Cahn-Hilliard, since the phase field dependent coefficients affects the whole system.
  The major difficulties cause from the phase field dependent coefficients
and from the coupling nonlinear terms.
 To the best of our knowledge, error estimates for fully discrete scheme of problem \eqref{eq2.1}-\eqref{eq4} are not yet set up \cite{Y2019}.
However, the analysis in \cite{2004Feng,Diegel2015,Caishen2018} cannot be easily developed to the fully
discrete form of  problem \eqref{eq2.1}-\eqref{eq4}, as the time-space discretization raises another difficulty, particularly
in deriving error estimates for the pressure.

In this paper, applying some useful techniques given
in references \cite{Tabata2005,Diegel2015,Caishen2018}, we analyze a fully discrete scheme for problem \eqref{eq2.1}-\eqref{eq4},
 which is comprised of conforming mixed finite element
method in space and the Euler semi-implicit discretization with a convex splitting method in time.
 The main purpose of this paper is to derive some optimal error estimates for the concentration field, the chemical potential,
 the velocity field, the magnetic field and the pressure for the Cahn-Hilliard-Magneto-hydrodynamics system in fully discrete form. Furthermore,
  the error estimates are established do not need any CFL conditions.
 The highlight of this
paper is to set up the error estimates for the fully discrete solution $(\phi^n_h,\mu^n_h,\textbf{u}^n_h,p^n_h,\textbf{B}^n_h)$ as follows:
\begin{align*}
\max_{1\leq n\leq N}\|\nabla(\phi(t_n)-\phi^n_h)\|+\Bigl(\Delta t\sum_{n=1}^{N}\|\nabla(\mu(t_n)-\mu^n_h)\|^2\Bigr)^{\frac{1}{2}}&\leq C\bigl(\Delta t+h^{k+1}\bigr),\\
\max_{1\leq n\leq N}\|\textbf{u}(t_n)-\textbf{u}^n_h\|+ \max_{1\leq n\leq N}\|\textbf{B}(t_n)-\textbf{B}^n_h\|&\leq C\bigl(\Delta t+h^{k+1}\bigr),\\
\Bigl(\Delta t\sum_{n=1}^{N}\|p(t_n)-p^n_h\|^2\Bigr)^{\frac{1}{2}}&\leq C\bigl(\Delta t+h^{k+1}\bigr).
\end{align*}
Numerical tests are given to
confirm the theoretical rates of convergence.

The paper is organized as follows.
In section 2 we present some preliminary results
and the well-posedness of weak solution
 of Cahn-Hilliard-Magneto-hydrodynamics system.  In section 3 we give
fully discrete scheme and obtain its unconditionally energy stability.
In section 4 we prove some optimal
error estimates for the concentration field, the chemical potential,
the velocity field, the magnetic field and the pressure. In section 5 we present
 some numerical tests to checked the theoretical results of the our scheme.

\section{\label{Sec2} Continuous problem}

\subsection{\label{Sec2.1} Preliminaries }
In this subsection we give some standard notations.
 Let $C^m(\Omega)$ ($m \in N$) be
the space of functions with up to $m$
times continuously differentiable in $\Omega$, and let $C^{m,1}(\Omega)$ be the space of
functions in $C^m(\Omega)$ that are Lipschitz continuous in $\Omega$.
 let $(L^p(\Omega),\|\cdot\|_{L^p})$ and $(W^{k,p}(\Omega),\|\cdot\|_{k,p})$ denote respectively the Lebesgue spaces and Sobolev spaces.
For simplicity, we denote $\|\cdot\|:=\|\cdot\|_{L^2}$,  and denote by $H^k(\Omega)$ the Sobolev space $W^{k,2}(\Omega)$.

To define a weak formulation of problem \eqref{eq2.1}-\eqref{eq4}, we introduce the following Sobolev spaces
\begin{align*}
&\cX:={H}_0^1(\Omega)^d=\bigl\{\textbf{v}\in H^1(\Omega)^d: \textbf{v}|_{\partial\Omega}=0\bigr\}, \\
&\cX_0:=\bigl\{\textbf{v}\in L^2(\Omega)^d: \nabla\cdot\textbf{v}=0,\textbf{v}|_{\partial\Omega}=0\bigr\}, \\
&\cW:={H}^1_n(\Omega)^d:=\bigl\{\textbf{w}\in H^1(\Omega)^d: \,\textbf{w}\cdot\textbf{n}|_{\partial\Omega}=0 \bigr\}, \\
&\cW_0:=\bigl\{\textbf{w}\in H^1(\Omega)^d: \,\curl\textbf{w}\times\textbf{n}|_{\partial\Omega}=0 \bigr\}, \\
&\cX_0:=\bigl\{\textbf{v}\in  \cX: \mbox{div}\textbf{v}=0\bigr\}, \\
&\cW_n: =\bigl\{\textbf{w}\in  \cW: \mbox{div}\textbf{w}=0 \bigr\},\\
&\cQ:=\bigl\{\varphi\in H^1(\Omega): \frac{\partial\varphi}{\partial\textbf{n}}|_{\partial\Omega}=0\bigr\}, \\
&\cM: =L_0^2(\Omega)=\bigl\{ q\in L^2(\Omega), \int_\Omega q d \textbf{x}=0 \bigr\}.
\end{align*}

 The following useful inequalities hold \cite{Temam1983,Girault1986,1991Gunzburger,2006Gerbeau,2007Feng,He2015}:
\begin{subequations}\label{eq6.1}
\begin{alignat}{2} \label{eq6.1a}
\|\textbf{v}\|&\leq C_{\Omega}\|\nabla\textbf{v}\|, \ \  &&\ \forall \ \textbf{v} \in \cX,\\\label{eq6.1b}
\|\textbf{v}\|_{L^6}&\leq   C_{\Omega}\|\nabla\textbf{v}\|, \ \ &&\ \forall \ \textbf{v} \in \cX,\\\label{eq6.1c}
\|\mathbb{D}\textbf{v}\|&\geq  c_0\|\nabla\textbf{v}\|, \ \ &&\ \forall \ \textbf{v} \in \cX,\\\label{eq6.1d}
\|\textbf{v}\|_{L^3}&\leq C_{\Omega}\|\textbf{v}\|^{\frac{6-d}{6}}\|\nabla\textbf{v}\|^{\frac{d}{6}}, \ \ &&\ \forall \ \textbf{v} \in \cX,\\\label{eq6.1e}
\|\textbf{v}\|_{L^4}&\leq  C_{\Omega}\|\textbf{v}\|^{\frac{4-d}{4}}\|\nabla\textbf{v}\|^{\frac{d}{4}}, \ \ &&\ \forall \ \textbf{v} \in \cX,\\\label{eq6.1f}
\|\textbf{v}\|_{L^\infty}&\leq   C_{\Omega}\|\textbf{v}\|_{H^1}^{1/2}\|\textbf{v}\|^{1/2}_{{H^2}},  &&\ \forall\ \textbf{v} \in {H}^2(\Omega)^d,\\\label{eq6.1g}
c_1 \|\nabla\textbf{B}\|^2&\leq  \|\curl \textbf{B}\|^2+\|\div \textbf{B}\|^2,&&\ \forall\ \textbf{B} \in \cW,\\\label{eq6.1h}
\|\curl \textbf{B}\|&\leq\sqrt{2}\|\nabla\textbf{B}\|,  \|\div\textbf{B}\|\leq \sqrt{d}\|\nabla\textbf{B}\|, &&\ \forall\ \textbf{B} \in {H}^1(\Omega)^d,\\\label{eq6.1i}
\|\phi\|_{L^p}&\leq   C_{\Omega}\|\phi\|_{H^1},  \quad (2\leq p\leq 6)\ \ &&\ \forall \ \textbf{v} \in \cQ,,\\\label{eq6.1j}
\|\phi\|_{L^3}&\leq   C_{\Omega}\|\phi\|^{\frac{6-d}{6}}\|\nabla\phi\|^{\frac{d}{6}}+C_{\Omega}\|\phi\|, \ \ &&\ \forall \ \textbf{v} \in \cQ,\\\label{eq6.1k}
\|\phi\|_{L^\infty}&\leq   C_{\Omega}\|\Delta\phi\|^{\frac{d}{2(6-d)}}\|\phi\|_{L^6}^{\frac{3(4-d)}{2(6-d)}}+C_{\Omega}\|\phi\|_{L^6}, \ \ &&\ \forall \ \textbf{v} \in \cQ,
\end{alignat}
\end{subequations}
where  $c_0$, $c_1$ and $C_{\Omega}$ are positive constants depending on $\Omega$.

 We define the following bilinear terms:
\begin{align*}
 a_\phi(\varphi;\phi,\theta)&=\int_\Omega \kappa(\varphi)\nabla\phi\cdot\nabla\theta d\textbf{x},
&&a_f(\varphi;\textbf{u},\textbf{v})=\int_\Omega2\nu(\varphi)\mathbb{D}(\textbf{u}):\mathbb{D}(\textbf{v})d\textbf{x}, \\
a_{B}(\varphi;\textbf{B},\textbf{H})&=\int_\Omega \eta(\varphi)\curl\textbf{B}\cdot\curl\textbf{H}d\textbf{x}
&&+\int_\Omega \eta(\varphi)\div\textbf{B}\cdot\div\textbf{H}d\textbf{x},\\
d(\textbf{v},q)&=\int_\Omega q\mbox{div}\textbf{v}d\textbf{x},
\end{align*}
and trilinear terms:
\begin{align*}
b(\textbf{w},\textbf{u},\textbf{v})&= \frac{1}{2}\int_\Omega[(\textbf{w}\cdot\nabla)\textbf{u}]\cdot \textbf{v}
 -[(\textbf{w}\cdot\nabla)\textbf{v}]\cdot \textbf{u} d\textbf{x}=\int_\Omega[(\textbf{w}\cdot\nabla)\textbf{u}]\cdot \textbf{v}
 +\frac{1}{2}[(\nabla\cdot\textbf{w})\textbf{u}]\cdot \textbf{v} d\textbf{x},\\
c_{\widehat{B}}(\textbf{H},\textbf{B},\textbf{v})&=\int_\Omega \textbf{H}\times \mbox{curl}\textbf{B}\cdot\textbf{v}d\textbf{x},\
\ \quad c_{\widetilde{B}}(\textbf{u},\textbf{B},\textbf{H}) =\int_\Omega(\textbf{u}\times \textbf{B})\cdot\curl\textbf{H}d\textbf{x}.
\end{align*}
In addition, using the definition of $b(\cdot,\cdot,\cdot)$, it follows that
\begin{align}\label{eq7}
b(\textbf{u},\textbf{v},\textbf{v})&=0,\  &&\textbf{u}\in \cX, \textbf{v}\in {H}^1(\Omega)^d.
\end{align}
Applying $(\textbf{B}\times \mbox{curl}\textbf{H},\textbf{v})=(\textbf{v}\times \textbf{H}, \mbox{curl}\textbf{B})$,  one finds that
\begin{align}\label{eq8}
c_{\widehat{B}}(\textbf{B},\textbf{B},\textbf{u})-c_{\widetilde{B}}(\textbf{u},\textbf{B},\textbf{B})=0,\ \ \textbf{u}\in \cX,  \textbf{B} \in \cW.
\end{align}
The bilinear term $d(\cdot,\cdot)$  satisfies the LBB condition \cite{Girault1986,Temam1983}:
\begin{align}\label{eq9}
\sup_{\textbf{v}\in \cX,\textbf{v}\neq \textbf{0}}\dfrac{d(\textbf{v},q)}{\|\textbf{v}\|_{1}}\geq \widehat{\beta}\|q\|,\ \ \  \forall q \in \cM,
\end{align}
where $\widehat{\beta}>0$ is constant depend on $\Omega$.

A weak formulation for \eqref{eq2.1}-\eqref{eq4} may be written as follows: find $(\phi,\mu,\textbf{u}, p,\textbf{B})$
 such that
 \begin{subequations}\label{eq10.1}
\begin{alignat}{2} \label{eq10.1a} (\partial_t \phi,\psi)+ \epsilon a_{\phi}(\phi;\mu,\psi)+(\nabla\phi\cdot\textbf{u},\psi)&=0, \\
\epsilon^{-1}(\phi^3-\phi,\theta)+\epsilon(\nabla\phi,\nabla\theta)&=(\mu,\theta), \\
(\partial_t\textbf{u},\textbf{v})+a_f(\phi;\textbf{u},\textbf{v})
+b(\textbf{u},\textbf{u},\textbf{v})+S_cc_{\widehat{B}}(\textbf{B},\textbf{B},\textbf{v})
\\\qquad\qquad\nonumber-d(\textbf{v},p)&=\lambda(\mu\nabla\phi,\textbf{v}), \\
d(\textbf{u},q)&=0, \\
(\partial_t\textbf{B},\textbf{H})+a_B(\phi;\textbf{B},\textbf{H})-c_{\widetilde{B}}(\textbf{u},\textbf{B},\textbf{H})&=0,
\end{alignat}
\end{subequations}
for all $(\psi,\theta, \textbf{v},q,\textbf{H})\in \cQ\times\cQ\times\cX\times \cM \times \cW$.

\subsection{\label{Sec2.3} Wellposedness of solution }

This subsection builds a well-posedness result of \eqref{eq10.1}.
 For simplicity, we consider the following problem: find $(\phi,\mu,\textbf{u},\textbf{B})\in  (\cQ, \cQ, \cX_0, \cW)$, such that
 \begin{subequations}\label{eq11.1}
\begin{alignat}{2} \label{eq11.1a} (\partial_t \phi,\psi)+ \epsilon a_{\phi}(\phi;\mu,\psi)+(\nabla\phi\cdot\textbf{u},\psi)&=0, \\
\epsilon^{-1}(\phi^3-\phi,\theta)+\epsilon(\nabla\phi,\nabla\theta)&=(\mu,\theta), \\
(\partial_t\textbf{u},\textbf{v})+a_f(\phi;\textbf{u},\textbf{v})
+b(\textbf{u},\textbf{u},\textbf{v})+S_cc_{\widehat{B}}(\textbf{B},\textbf{B},\textbf{v})&=\lambda(\mu\nabla\phi,\textbf{v}), \\
(\partial_t\textbf{B},\textbf{H})+a_B(\phi;\textbf{B},\textbf{H})-c_{\widetilde{B}}(\textbf{u},\textbf{B},\textbf{H})&=0,
\end{alignat}
\end{subequations}
for all $(\psi,\theta, \textbf{v},\textbf{H})\in \cQ\times\cQ\times\cX_0\times\cW$.

The nest theorems give a well-posedness result for weak solution of problem \eqref{eq11.1}.
They were proved by the similar lines as in \cite{Lorca1999,Bermudez2010,2014Han}.
 Thus, we skip the proofs of the following theorems.
\begin{theorem} Suppose that the initial conditions $\phi_0, \textbf{u}_0, \textbf{B}_0$ satisfy
\begin{align}\label{eq12}
\phi_0\in H^{1}(\Omega), \qquad \textbf{u}_0\in L^2(\Omega)^d, \qquad \textbf{B}_0\in L^2(\Omega)^d.
 \end{align}
Furthermore, assume that the given functions $\kappa$, $\nu$ and $\eta$ satisfy
 \begin{align}\label{eq13}
 \kappa, \nu, \eta \in C(\overline{\Omega}\times [0,T]\times R; R^{+}).
 \end{align}
 Then problem \eqref{eq11.1} has at least one solution
 $(\phi,\mu,\textbf{u},\textbf{B})$ such that
  \begin{align}\label{eq14}
 \phi&\in L^{\infty}(0,T,H^{1}(\Omega))\cap L^2(0,T,H^{3}), \\\nonumber
 \textbf{u}&\in L^{\infty}(0,T,L^2(\Omega)^d)\cap L^2(0,T,\cX), \\\nonumber
 \textbf{B}&\in L^{\infty}(0,T,L^2(\Omega)^d)\cap L^2(0,T,\cW)\\\nonumber
  \mu&\in L^2(0,T,\cQ).
 \end{align}
\end{theorem}
\begin{theorem}  Taking the place of \eqref{eq13} by
 \begin{align}\label{eq15}
 \nu, \eta, \kappa \in C^{0,1}(\overline{\Omega}\times [0,T]\times R; R^{+}),
 \end{align}
and retaining the other assumptions in \textbf{Theorem 1}. In addition,  we suppose $(\phi,\mu,\textbf{u},\textbf{B})$ satisfy
 \begin{align}\label{eq16}
 \phi\in L^2(0,T,W^{1,\infty}),\ \quad
\textbf{u} \in L^2(0,T,W^{1,\infty}(\Omega)^d),\,\quad
\textbf{B}\in L^2(0,T,W^{1,\infty}(\Omega)^d).
 \end{align}
 Then problem \eqref{eq11.1} has a unique solution $(\phi,\mu,\textbf{u},\textbf{B})$.
\end{theorem}
\section{\label{Sec3} Fully discrete scheme }

In this section, we give a fully discrete scheme
based on applying conforming finite element method
in space and Euler semi-implicit discretization with a convex splitting method in time
for \eqref{eq10.1} and obtain some unconditionally energy stability.

Let $\textbf{{K}}_h$ be a conforming, quasi-uniform family of
triangulations of $\Omega$ into triangles
when $d$=2 and tetrahedra when $d=3$, respectively.
Furthermore, we introduce four finite element spaces
$\cX_h, \cM_h, \cW_h, \cQ_h$  with  $\cX_h \subset \cX, \cM_h \subset \cM,
 \cW_h \subset \cW, \cQ_h \subset \cQ$  as follows.
 \begin{align*}
\cX_h&=\bigl\{\textbf{v}\in \textbf{C}^0(\overline{\Omega})^d\cap \cX:\textbf{v}|_{\mathbb{K}}\in {P}_{r+1}(\mathbb{K})^d,\forall\, \textbf{{K}}\in \textbf{{K}}_h\bigr\},\\
 \cM_h&=\bigl\{q\in C^0(\overline{\Omega})\cap\cM:q|_{\mathbb{K}}\in P_r(\mathbb{K}),\forall\, \textbf{K}\in \textbf{{K}}_h\bigr\}, \\
\cW_h&=\bigl\{\textbf{w}\in \textbf{C}^0(\overline{\Omega})^d\cap \cW:\textbf{w}|_{\mathbb{K}}\in {P}_{r+1}(\mathbb{K})^d,\forall\, \textbf{{K}}\in \textbf{{K}}_h\bigr\}, \\
 \cQ_h&=\bigl\{\omega\in C^0(\overline{\Omega})\cap\cQ:\omega|_{\mathbb{K}}\in P_{r+1}(\mathbb{K}),\forall\, \textbf{{K}}\in \textbf{{K}}_h\bigr\}\\
\cX_{0h}&=\bigl\{\textbf{v}\in  \cX_h:d(\textbf{v},w)=0,\forall\, w\in \cM_h\bigr\}.
 \end{align*}

As is noted that the $(\cX_h,\cM_h)$ is Taylor-Hood finite element pair.
Therefore, the finite element pair $(\cX_h, \cM_h)$ satisfies the discrete LBB condition \cite{Temam1983,Girault1986}:

\textbf{Assumption A1:} The following discrete LBB condition holds:
\begin{align*}
\exists {\beta}_0>0, \quad \sup_{\textbf{v}_h\in\cX_h,\textbf{v}_h\neq 0}\dfrac{d(\textbf{v}_h,q_h)}{\|\nabla\textbf{v}_h\|}\geq {\beta}_0\|q_h\|, \ \forall q_h \in {\cM}_h.
 \end{align*}

 Moreover, we suppose the finite element spaces satisfy following inverse inequality and finite element approximation properties:

 \textbf{Assumption A2:}  The following inverse inequality holds
\begin{align*}
 \|{v}_h\|_{m,q}&\leq  Ch^{l-m+d(\frac1q-\frac1p)}\|{v}_h\|_{l,p},  \forall \ {v}_h \in \cX_h,\cW_h\, \mbox{or}\,\cQ_h,\\&
 \quad  0 \leq l \leq m \leq 1, \quad 1 \leq p\leq q \leq \infty.
 \end{align*}

 \textbf{Assumption A3:} There exists $k\geq 1$, such that for all $1\leq l\leq k$,
\begin{align*}
  \inf_{\textbf{v}_h \in\cX_h}\bigl[\|\textbf{v}-\textbf{v}_h\|+h\|\nabla(\textbf{v}-\textbf{v}_h)\|\bigr]&\leq Ch^{l+1}\|\textbf{v}\|_{l+1}, \ &&\forall \textbf{v} \in H^{l+1}(\Omega)^d,\\
   \inf_{\textbf{B}_h\in\cW_h}\bigl[\|\textbf{B}-\textbf{B}_h\|+h\|\nabla(\textbf{B}-\textbf{B}_h)\|\bigr]&\leq Ch^{l+1}\|\textbf{B}\|_{l+1}, \ &&\forall \textbf{B} \in H^{l+1}(\Omega)^d,\\
    \inf_{\varphi_h\in\cQ_h}\bigl[\|\varphi-\varphi_h\|+h\|\nabla(\varphi-\varphi_h)\|\bigr]&\leq Ch^{l+1}\|\varphi\|_{l+1}, \  &&\forall \varphi \in H^{l+1}(\Omega),\\
     \inf_{q_h\in\cM_h}\|q-q_h\| &\leq Ch^{l}\|q\|_{l}, \ && \forall q \in H^{l}(\Omega).
 \end{align*}

 Let $N$ be a positive integer and $0=t_0<t_1<\cdots<t_N=T$ be a uniform partition of $[0,T]$, with $\Delta t=t_i-t_{i-1}$ for $i=1,\cdots,N$.
Let us denote by $\phi^n$ the valve $\phi(t_n)$ at the time $t_n$, and denote $\delta_t{\phi}:=\frac{\phi^n-\phi^{n-1}}{\Delta t}$.
Denote $\kappa^n:=\kappa(\textbf{x},t_n,\phi)$ $\nu^n:=\nu(\textbf{x},t_n,\phi)$ and $\eta^n:=\eta(\textbf{x},t_n,\phi)$.

 With above preparation, we give the fully discrete scheme of problem \eqref{eq10.1}.

\textbf{Scheme 3.1:}  Given $(\phi^0_h,\textbf{u}^0_h,\textbf{B}^0_h) \in \cQ_h\times\cX_h\times \cW_h$,
find $(\phi^n_h,\mu^n_h,\textbf{u}^n_h, p^n_h, \textbf{B}^n_h) \in \cQ_h\times \cQ_h\times{\cX}_h\times \cM_h\times \cW_h$, such that
 \begin{subequations}\label{eq17.1}
\begin{alignat}{2} \label{eq17.1a}
(\delta_t \phi^n_h,\psi_h)+ \epsilon a_{\phi}(\phi^{n-1}_h;\mu^{n}_h,\psi_h)+(\nabla\phi^{n-1}_h\cdot\textbf{u}^n_h,\psi_h)&=0, \\\label{eq17.1b}
 \epsilon^{-1}\bigl((\phi_h^{n})^3-\phi^{n-1}_h,\theta_h\bigr)+\epsilon(\nabla\phi_h^n,\nabla\theta_h)&=(\mu^n_h,\theta_h), \\\label{eq17.1c}
(\delta_t\textbf{u}^n_h,\textbf{v}_h)+a_f(\nu^n(\phi^{n-1}_h);\textbf{u}^n_h,\textbf{v}_h)
+b_f(\textbf{u}^{n-1}_h,\textbf{u}^n_h,\textbf{v}_h)\\\nonumber
+S_cb_{\widehat{B}}(\textbf{B}^{n-1}_h,\textbf{B}^n_h,\textbf{v}_h)-d(\textbf{v}_h,p^n_h)&=\lambda(\mu_h^n\nabla\phi_h^{n-1},\textbf{v}_h), \\\label{eq17.1d}
d(\textbf{u}^n_h,q_h)&=0,\\\label{eq17.1e}
(\delta_t\textbf{B}^n_h,\textbf{H}_h)+a_B(\eta^n(\phi^{n-1}_h);\textbf{B}^{n}_h,\textbf{H}_h)
-c_{\widetilde{B}}(\textbf{u}^{n}_h,\textbf{B}^{n-1}_h,\textbf{H}_h)&=0,
\end{alignat}
\end{subequations}
for all $(\psi_h,\theta_h,\textbf{v}_h,q_h,\textbf{H}_h)\in \cQ_h\times \cQ_h\times{\cX}_h\times \cM_h \times \cW_h$.

 In addition, the given functions  $\kappa, \nu, \eta$  satisfy following assumption.

  \textbf{Assumption A4:} For all $(\textbf{x},t,\phi) \in \overline{\Omega}\times (0,T] \times R$, the given functions $\kappa, \nu, \eta$  satisfy
  \begin{align*}
   0<\kappa_1\leq \kappa(\textbf{x},t,\phi)&\leq \kappa_2,\\
  0<\nu_1\leq \nu(\textbf{x},t,\phi)&\leq \nu_2,\\
 0<\eta_1\leq \eta(\textbf{x},t,\phi)&\leq \eta_2.
 \end{align*}

 Denote $\dot{S}_h=\cQ_h\cap\cM_h$. Let $\hat{\kappa} \in C^{0,1}(\overline{\Omega};R^{+})$ satisfy \textbf{Assumption A4}, we
define the operator  $T_h:\dot{S}_h\rightarrow \dot{S}_h$ by the following weak problem: given $\sigma\in \dot{S}_h$, find $T_h(\sigma) \in \dot{S}_h$
such that
 \begin{align}\label{eq18}
a_{\phi}(\hat{\kappa};T_h(\sigma),\varrho)=(\sigma,\varrho), \ \ \quad \forall \varrho\in \dot{S}_h.
\end{align}
As in the work \cite{Diegel2015}, we have similar results.
\begin{lemma} Let $\zeta,\vartheta\in \dot{S}_h$ and denote
 \begin{align*}
 (\zeta,\vartheta)_{-1,h}:= a_{\phi}(\hat{\kappa};T_h(\zeta),T_h(\theta))=(\zeta,T_h(\theta))=(T_h(\zeta),\theta),
 \end{align*}
 and define the following negative norm:
  \begin{align*}
\|\zeta\|_{-1,h}:=\sqrt{(\zeta,\zeta)_{-1,h}}=\sup_{0\neq \theta \in \dot{S}_h}\frac{(\zeta,\theta)}{\|\nabla\theta \|}.
 \end{align*}
 Therefore, $\forall \, \vartheta \in \cX_h$ and $\forall\, \zeta \in \dot{S}_h$, we have
   \begin{align*}
  |(\zeta,\vartheta)|\leq \|\zeta\|_{-1,h}\|\nabla\vartheta\|.
 \end{align*}
 Then the following Poincar\'{e} inequality holds:
   \begin{align*}
 \|\vartheta\|_{-1,h}\leq C\|\vartheta\|, \ \ \qquad \forall\, \vartheta \in \dot{S}_h.
 \end{align*}
\end{lemma}

We now are ready to state and prove unconditionally energy stability for the fully discrete scheme \eqref{eq17.1}.
\begin{theorem}  Suppose \textbf{Assumptions A1-A4} hold. Then the finite element approximate
 solution $(\phi^n_h,\mu^n_h,\textbf{u}^n_h,p^n_h,\textbf{B}^n_h)$  of problem
 \eqref{eq17.1} satisfy the following discrete energy law:
  \begin{align*}
\textbf{E}(\phi^{n}_h,{\textbf{u}}^{n}_h,{\textbf{B}}^{n}_h)
&+\Delta t\sum^n_{m=1}\bigl[\kappa_1\epsilon\|\nabla{\mu}^m_h\|^2+\Delta tc_0\nu_1\|\nabla{\textbf{u}}^n_h\|^2+S_c\Delta tc_1\eta_1\|\nabla{\textbf{B}}^n_h\|^2\bigl]\\&
+\Delta t^2\sum^n_{m=1}\bigl[\frac{\epsilon}{2}\|\nabla\delta_t\phi^n_h\|^2+\frac{1}{2\lambda}\|\delta_t{\textbf{u}}^n_h\|^2+\frac{S}{2}\|\delta_t{\textbf{B}}^n_h\|^2
\\&+\frac{1}{4\epsilon}\|\delta_t(\phi^n_h)^2\|^2+\frac{1}{2\epsilon}\|\phi^n_h\delta_t\phi^n_h\|^2+\frac{1}{2\epsilon}\|\delta_t\phi^n_h\|^2\bigr]\leq\textbf{E}(\phi^{0}_h,{\textbf{u}}^{0}_h,{\textbf{B}}^{0}_h).
 \end{align*}
\end{theorem}
\noindent \textit{Proof:}\quad Taking $\psi_h=\mu^n_h$ in \eqref{eq17.1a}, $\theta_h=\delta_t\phi^n_h$ in \eqref{eq17.1b}, $\textbf{v}_h=\textbf{u}^n_h$ in \eqref{eq17.1c}, $q_h=p^n_h$ in \eqref{eq17.1d}, $\textbf{H}_h=\textbf{B}^n_h$ in \eqref{eq17.1e}, one finds that
\begin{align*}
(\delta_t\phi^n_h,\mu^n_h)+ \epsilon a_{\phi}(\phi^{n-1}_h;\mu^{n}_h,\mu^n_h)+(\nabla\phi^{n-1}_h\cdot\textbf{u}^n_h,\mu^n_h)&=0, \\
 \epsilon^{-1}\bigl((\phi_h^{n})^3-\phi^{n-1}_h,\delta_t\phi^n_h\bigr)+\epsilon(\nabla\phi_h^n,\nabla\delta_t\phi^n_h)&=(\mu^n_h,\delta_t\phi^n_h), \\
(\delta_t\textbf{u}^n_h,\textbf{u}^n_h)+a_f(\nu^n(\phi^{n-1}_h);\textbf{u}^n_h,\textbf{u}^n_h)
+b_f(\textbf{u}^{n-1}_h,\textbf{u}^n_h,\textbf{u}^n_h)\\\nonumber
+S_cb_{\widehat{B}}(\textbf{B}^{n-1}_h,\textbf{B}^n_h,\textbf{u}^n_h)-d(\textbf{u}^n_h,p^n_h)&=\lambda(\mu_h^n\nabla\phi_h^{n-1},\textbf{u}^n_h), \\
d(\textbf{u}^n_h,p_h)&=0,\\
(\delta_t\textbf{B}^n_h,\textbf{B}^n_h)+a_B(\eta^n(\phi^{n-1}_h);\textbf{B}^{n}_h,\textbf{B}^n_h)
-c_{\widetilde{B}}(\textbf{u}^{n}_h,\textbf{B}^{n-1}_h,\textbf{B}^n_h)&=0.
 \end{align*}
 Using the elementary identity
 \begin{align*}
 2a(a-b)=a^2-b^2+(a-b)^2,\ \ \ \mbox{for all} \ a, b, \in \R^d.
 \end{align*}
It follows that
\begin{align*}
\epsilon(\nabla\phi_h^n,\nabla\delta_t\phi^n_h)&=\frac{\epsilon}{2}\bigl[\delta_t\|\nabla\phi^n_h\|^2+\Delta t\|\nabla\delta_t\phi^n_h\|^2\bigr],\\
 \epsilon^{-1}\bigl((\phi_h^{n})^3-\phi^{n-1}_h,\delta_t\phi^n_h\bigr)&=\frac{\epsilon^{-1}}{4}\delta_t\|(\phi^n_h)^2-1\|^2\\
 &\quad+\frac{\epsilon^{-1}\Delta t}{4}\bigl[\|\delta_t(\phi^n_h)^2\|^2+2\|\phi^n_h\delta_t\phi^n_h\|^2+2\|\delta_t\phi^n_h\|^2\bigr],\\
(\delta_t\textbf{u}^n_h,\textbf{u}^n_h)&=\frac{1}{2}\bigl[\delta_t\|\textbf{u}^n_h\|^2+\Delta t\|\delta_t\textbf{u}^n_h\|^2\bigr],\\
(\delta_t\textbf{B}^n_h,\textbf{B}^n_h)&=\frac{1}{2}\bigl[\delta_t\|\textbf{B}^n_h\|^2+\Delta t\|\delta_t\textbf{B}^n_h\|^2\bigr].
 \end{align*}
Then using the operator $\Delta t\sum^m_{n=1}$ to the combined equations, the desired result is derived.
 The proof is completed.
 $$\eqno\Box$$

The following Theorem give some bounds for the fully discrete solution. It is very important to derive optimal error estimates of problem \eqref{eq17.1}.
 \begin{theorem} \cite{Diegel2015,He2015,Y2019} Suppose \textbf{Assumptions A1-A4} hold.
 Then the fully discrete solution $(\phi^n_h,\mu^n_h,\textbf{u}^n_h,p^n_h,\textbf{B}^n_h)$  of \eqref{eq17.1} satisfy the following bounds:
  \begin{align}\label{eq19}
& \max_{0\leq m\leq n}\bigl[\|\textbf{u}^n_h\|^2+\|\textbf{B}^n_h\|^2+\|\nabla\phi^n_h\|^2+\|(\phi^n_h)^2-1\|^2+\|\phi^n_h\|^4 \\\nonumber&\qquad\qquad\qquad\qquad
+\|\phi^n_h\|^2+\|\phi^n_h\|^2_{H^1}\bigr]\leq C,\\\label{eq20}
 &\Delta t\sum_{m=1}^n\bigl[\|\nabla\textbf{u}^n_h\|^2+\|\nabla\textbf{B}^n_h\|^2+\|\nabla\mu^n_h\|^2\bigr]\leq C,\\\label{eq21}
  &\sum_{m=1}^n\bigl[\|\nabla(\phi^n_h-\phi^{n-1}_h)\|^2+\|\phi^n_h-\phi^{n-1}_h\|^2 +\|\phi^n_h(\phi^n_h-\phi^{n-1}_h)\|^2  \\\nonumber&\qquad\qquad\qquad\qquad
+\|(\phi^n_h)^2-(\phi^{n-1}_h)^2\|^2+\|\phi^n_h-\phi^{n-1}_h\|_{-1,h}^2 \\\nonumber&\qquad\qquad\qquad\qquad
+\|\textbf{u}^n_h-\textbf{u}^{n-1}_h\|^2+\|\textbf{B}^n_h-\textbf{B}^{n-1}_h\|^2\bigr]\leq C,\\\label{eq22}
   &\Delta t\sum_{m=1}^n\bigl[\|\delta_t\phi^n_h\|_{H^{-1}}^2+\|\delta_t\phi^n_h\|_{-1,h}^2+\|\nabla_h\phi^n_h\|^2+\|\mu^n_h\|^2
   \\\nonumber&\qquad\qquad\qquad\qquad+\|\mu^n_h\|^{\frac{4(6-d)}{d}}_{L^{\infty}}+\|\delta_t\phi^n_h\|^2\bigr]\leq C,\\\label{eq23}
    & \max_{0\leq m\leq n}\bigl[\|\mu^n_h\|^2+\|\nabla_h\phi^n_h\|^2+\|\mu^n_h\|^{\frac{4(6-d)}{d}}_{L^{\infty}}\bigr]\leq C,\\\label{eq24}
   &\Delta t^2\sum_{m=1}^n\bigl[\|\delta_t\textbf{u}^n_h\|^2+\|\delta_t\textbf{B}^n_h\|^2\bigr]\leq C.
 \end{align}
\end{theorem}


\section{\label{Sec4}Error analysis}

In this section we drive some optimal error estimates of the fully discrete scheme \eqref{eq17.1}.
Therefore, we assume that the solution $(\phi,\mu,\textbf{u},p,\textbf{B})$ has the following regularity.

 \textbf{Assumption A5:} The weak solution $(\phi,\mu,\textbf{u},p,\textbf{B})$ of
 Cahn-Hilliard-Magneto-hydrodynamics \eqref{eq2.1}-\eqref{eq4} is sufficiently smooth such that
 \begin{align*}
 &\partial_{tt}\phi \in L^{\infty}(0,T;H^1(\Omega)), \quad&&\phi, \partial_{t}\phi\in L^{\infty}(0,T;H^{k+2}(\Omega)\cap W^{1,\infty}(\Omega)),\\
 &\partial_{tt}\textbf{u} \in L^{\infty}(0,T;H^1(\Omega)^d),\quad &&\textbf{u}, \partial_{t}\textbf{u}\in L^{\infty}(0,T;H^{k+2}(\Omega)^d\cap W^{1,\infty}(\Omega)^d),\\
&\partial_{tt}\textbf{B} \in L^{\infty}(0,T;H^1(\Omega)^d), \quad &&\textbf{B}, \partial_{t}\textbf{B}\in L^{\infty}(0,T;H^{k+2}(\Omega)^d\cap W^{1,\infty}(\Omega)^d),\\
& \mu \in  L^{\infty}(0,T;H^{k+2}(\Omega)),  \quad &&p,\partial_tp\in L^{\infty}(0,T;H^{k+1}(\Omega)).
 \end{align*}

We need the following technical lemmas to obtain the rate of convergence of the fully discrete scheme \eqref{eq17.1}.
\begin{lemma}\cite{Diegel2015} Assume $\varpi \in H^1(\Omega)$, and $\upsilon\in \dot{S}_h$, Then there exists $C>0$,
independent of $h$, such that
 \begin{align}\label{eq25}
|(\varpi,\upsilon)|\leq C\|\nabla\varpi\|\|\upsilon\|_{-1,h}.
 \end{align}
\end{lemma}
\begin{lemma}\cite{Diegel2015} Assume that $(\phi,\mu,\textbf{u},p,\textbf{B})$, $(\phi_h^n,\mu_h^n,\textbf{u}_h^n,p_h^n,\textbf{B}_h^n)$
are weak solution to problem \eqref{eq10.1} and \eqref{eq17.1}, respectively. Then for any $h$,
$\Delta t>0$, the following inequality holds
\begin{align}\label{eq26}
\|\nabla(\phi^3-(\phi_h^n)^3)\|\leq C\|\nabla(\phi-{\phi_h^n})\|.
 \end{align}
\end{lemma}

  Let $\hat{\nu},\hat{\eta},\hat{\kappa} \in C^{0,1}(\overline{\Omega};R^{+})$ satisfy \textbf{Assumption A4}, we introduce the following four projections.
  Let $(\textbf{u},p)\in \cX\times \cM$, $\textbf{B} \in \cW$, $\theta, \phi \in \cQ$,
  we introduce Stokes projection $(R_h\textbf{u},J_hp) \in \cX_{h}\times \cM_h$ as the solution of the weak problem as follows
\begin{align}\label{eq27}
a_f(\hat{\nu};\textbf{u}-R_h\textbf{u},\textbf{v})+d(\textbf{v},p-J_hp)&=0, \ \ \forall\ \textbf{v}\in \cX_h,\\\label{eq28}
d(\textbf{u}-R_h\textbf{u},q)&=0, \ \ \forall\ q\in \cM_h.
 \end{align}
We have the following approximation result
\begin{align}\label{eq29}
\|\textbf{u}-R_h\textbf{u}\|+h\bigl(\|\nabla(\textbf{u}-R_h\textbf{u})\|+\|p-J_hp\|\bigr)&\leq Ch^{k+2}\bigl(\|\textbf{u}\|_{k+2}+\|p\|_{k+1}\bigr).
 \end{align}
Similarly, Maxwell projection $R_{mh}\textbf{B}\in \cW_{h}$ satisfying
\begin{align}\label{eq30}
a_B(\hat{\eta};\textbf{B}-R_{mh}\textbf{B},\textbf{H})&=0, \ &&\forall \ \textbf{H}\in \cW_h,\\\label{eq31}
\|\textbf{B}-R_{mh}\textbf{B}\|+h\|\nabla(\textbf{B}-R_{mh}\textbf{B})\|&\leq Ch^{k+2}\|\textbf{B}\|_{k+2},
 \end{align}
Ritz projections $r_{h}\theta \in \cQ_h$ and $\widehat{r}_{h}\theta \in \cQ_h$  satisfying
\begin{align}\label{eq32}
a_{\phi}(\hat{\kappa};\theta-r_h\theta,\psi)&=0, \ &&\forall \ {\psi}\in \cQ_h,\\\label{eq33}
\|\theta-r_h\theta\|+h\|\nabla(\theta-r_h\theta)\|&\leq Ch^{k+2}\|\theta\|_{k+2},\\\label{eq34}
(\nabla(\phi-\widehat{r}_h\phi),\nabla\psi)&=0, \ &&\forall \ {\psi}\in \cQ_h,\\\label{eq35}
\|\phi-\widehat{r}_h\phi\|+h\|\nabla(\phi-\widehat{r}_h\phi)\|&\leq Ch^{k+2}\|\phi\|_{k+2}.
 \end{align}

For convenience, let us define, for every $n\geq 0$,
\begin{align*}
e^n_\phi :&= \phi^n - \phi^n_h,\ &&e^n_\mu := \mu^n -\mu^n_h,\\
e^n_\textbf{u} :&= \textbf{u}^n - \textbf{u}^n_h,\ &&e^n_p := p^n -p^n_h,\\
e^n_\textbf{B} :&= \textbf{B}^n- \textbf{B}^n_h.
 \end{align*}
Which gives
\begin{align*}
 e^n_\phi&=\eta^n_\phi-\epsilon^n_\phi,  &&e^n_\mu=\eta^n_\mu-\epsilon^n_\mu,\\
 e^n_\textbf{u}&=\eta^n_\textbf{u}-\epsilon^n_\textbf{u},  &&e^n_p=\eta^n_p-\epsilon^n_p,\\
  e^n_\textbf{B}&=\eta^n_\textbf{B}-\epsilon^n_\textbf{B},
 \end{align*}
 with
 \begin{align*}
 &\eta^n_\phi:=\phi^n-r_h\phi^n \in \cQ, &&\epsilon_\phi^n:=\phi^n_h-r_h\phi^n\in \cQ_h,\\
  &\eta^n_\mu:=\mu^n-\widehat{r}_h\mu^n \in \cQ, &&\epsilon_\mu^n:=\mu^n_h-\widehat{r}_h\mu^n\in \cQ_h,\\
 &\eta^n_\textbf{u}:=\textbf{u}^n-R_h\textbf{u}^n \in \cX,\ \ \    && \epsilon_\textbf{u}^n:=\textbf{u}_h^n-R_h\textbf{u}^n \in \cX_h,\\
&\eta^n_p:=p^n-J_hp^n \in \cM,                     &&\epsilon_p^n:=p_h^n-J_hp^n \in \cM_h,\\
 &\eta^n_\textbf{B}:=\textbf{B}^n-R_{mh}\textbf{B}^n \in \cW,           &&\epsilon_\textbf{B}^n:= \textbf{B}^n_h-R_{mh}\textbf{B}^n \in \cW_h.
 \end{align*}

We are now in a position to give and derive the first main theorem of this section
for the concentration field, the chemical potential, the velocity field and the magnetic field.

\begin{theorem} Suppose that \textbf{Assumption A1-A5} hold with two positive constants $h_0,\Delta t_0$,
and suppose the scheme \eqref{eq17.1} is initialized such that
 \begin{align}\label{eq36}
\|\phi^0-\phi^0_h\|+\|\textbf{u}^0-\textbf{u}^0_h\|+ \|\textbf{B}^0-\textbf{B}^0_h\|\leq C\bigl(\Delta t+h^{k+1}\bigr).
 \end{align}
 For $h\in(0,h_0]$ and $\Delta t\in(0,\Delta t_0]$, the finite element approximate solution $(\phi^n_h,\mu^n_h,\textbf{u}^n_h,p^n_h,\textbf{B}^n_h)$ of \eqref{eq17.1} satisfy the following error estimates:
\begin{align*}
\max_{1\leq n\leq N}\|\nabla (\phi(t_n)-\phi^n_h)\|^2
+\Delta t\sum_{n=1}^{N}\|\nabla (\mu(t_n)-\mu^n_h)\|^2 &\leq C\bigl(\Delta t^2+h^{2k+2}\bigr),\\
\max_{1\leq n\leq N}\|\textbf{u}(t_n)-\textbf{u}^n_h\|^2
+\Delta t\sum_{n=1}^{N}\|\nabla (\textbf{u}(t_n)-\textbf{u}^n_h)\|^2 &\leq C\bigl(\Delta t^2+h^{2k+2}\bigr),\\
\max_{1\leq n\leq N}\|\textbf{B}(t_n)-\textbf{B}^n_h\|^2
+\Delta t\sum_{n=1}^{N}\|\nabla (\textbf{B}(t_n)-\textbf{B}^n_h)\|^2 &\leq C\bigl(\Delta t^2+h^{2k+2}\bigr).
 \end{align*}
\end{theorem}
\noindent \textit{Proof:}\quad Appying\eqref{eq10.1}, \eqref{eq17.1}, \eqref{eq27}-\eqref{eq28}, \eqref{eq30}, \eqref{eq32} and \eqref{eq34},  we can obtain the following error equations:
 \begin{subequations}\label{eq37.1}
\begin{alignat}{2} \label{eq37.1a}
(\delta_t\epsilon_\phi^n,\psi_h)+\epsilon a_\phi(\kappa^n(\phi^{n-1}_h);\epsilon_\mu^n,\psi_h)&=({\Lambda}_h^n,\psi_h)+(R^n_\phi,\psi_h),\\\label{eq37.1b}
\epsilon(\nabla\epsilon_\phi^n,\nabla\theta_h)-(\epsilon_\mu^n,\theta_h)&=(\widehat{\Lambda}_h^n,\theta_h),\\\label{eq37.1c}
 (\delta_t\epsilon_\textbf{u}^n,\textbf{v}_h)+a_f(\nu^n(\phi^{n-1}_h);\epsilon_\textbf{u}^n,\textbf{v}_h)
-d(\textbf{v}_h,\epsilon_p^n)&=(\Phi_h^n,\textbf{v}_h)+(R^n_\textbf{u},\textbf{v}_h),\\\label{eq37.1d}
d(\epsilon_\textbf{u}^n,q_h)&=0,\\\label{eq37.1e}
S_c(\delta_t\epsilon_\textbf{B}^n,\textbf{H}_h)+S_ca_B(\eta^n(\phi^{n-1}_h);\epsilon_\textbf{B}^n,\textbf{H}_h)
&=(\widehat{\Phi}_h^n,\textbf{H}_h)+S_c(R^n_\textbf{B},\textbf{H}_h).
\end{alignat}
\end{subequations}
where $(R^n_\phi,\psi_h)$, $(R^n_\textbf{u},\textbf{v}_h)$, $(R^n_\textbf{B},\textbf{H}_h)$,
$({\Lambda}_h^n,\psi_h)$, $(\widehat{\Lambda}_h^n,\theta_h)$, $(\Phi_h^n,\textbf{v}_h)$ and $(\widehat{\Phi}_h^n,\textbf{H}_h)$  are denoted by
 \begin{align}\label{eq38}
(R^n_\phi,\psi_h)&:=\bigl(\partial_t\phi^n-\delta_tr_h{\phi}^n,\psi_h\bigr),\\\label{eq39}
(R^n_\textbf{u},\textbf{v}_h)&:=\bigl(\partial_t\textbf{u}^n-\delta_tR_h{\textbf{u}}^n,\textbf{v}_h\bigr),\\\label{eq40}
(R^n_\textbf{B},\textbf{H}_h)&:=\bigl(\partial_t\textbf{B}^n-\delta_tR_{mh}{\textbf{B}}^n,\textbf{H}_h\bigr),\\\label{eq41}
({\Lambda}_h^n,\psi_h)&:=a_\phi(\kappa^n(\phi^n); \mu^n,\psi_h)-a_\phi(\kappa^n(\phi^{n-1}_h); \mu^n, \psi_h)\\&\nonumber
  \quad+(\textbf{u}^n\cdot\nabla \phi^n,\psi_h)-(\textbf{u}^{n}_h\cdot\nabla \phi^{n-1}_h,\psi_h),\\\label{eq42}
(\widetilde{\Lambda}_h^n,\theta_h)&:=(\eta_\mu^n,\theta_h)-\frac{\Delta t}{\epsilon}(\delta_t\phi^n,\theta_h)
\\&\nonumber\quad-\epsilon^{-1}(\phi^{n-1}-\phi^{n-1}_h,\theta_h)+\epsilon^{-1}\bigl((\phi^{n})^3-(\phi^{n}_h)^3,\theta_h\bigr),\\\label{eq43}
(\Phi_h^n,\textbf{v}_h)&: =a_f(\nu^n(\phi^n); \textbf{u}^n, \textbf{v}_h)-a_f(\nu^n(\phi^{n-1}_h); \textbf{u}^n, \textbf{v}_h)\\&\nonumber\quad
 +b(\textbf{u}^n,\textbf{u}^n,\textbf{v}_h)-b(\textbf{u}^{n-1}_h,\textbf{u}^n_h,\textbf{v}_h)\\&\nonumber
  \quad+S_cc_{\widehat{B}}(\textbf{B}^n, \textbf{B}^n,\textbf{v}_h)-S_cc_{\widehat{B}}(\textbf{B}^{n-1}_h, \textbf{B}^n_h,\textbf{v}_h)\\&\nonumber
  \quad+\lambda(\mu_h^n\nabla\phi_h^{n-1},\textbf{v}_h)-\lambda(\mu^n\nabla\phi^{n},\textbf{v}_h),\\\label{eq44}
(\widehat{\Phi}_h^n,\textbf{H}_h)&:=S_ca_B(\eta^n(\phi^n); \textbf{B}^n, \textbf{H}_h)-S_ca_B(\eta^n(\phi^{n-1}_h); \textbf{B}^n, \textbf{H}_h)\\&\nonumber
  \quad+S_cc_{\widehat{B}}(\textbf{u}^{n}_h, \textbf{B}^{n-1}_h,\textbf{H}_h)-S_cc_{\widehat{B}}(\textbf{u}^n, \textbf{B}^n,\textbf{H}_h).
 \end{align}
Setting $\psi_h=\epsilon_\phi^n \in \cQ_h$,
$\theta_h=\kappa_1\epsilon_\mu^n \in \cQ_h$
in \eqref{eq37.1a}-\eqref{eq37.1b}, respectively,  we have
 \begin{subequations}\label{eq45.1}
\begin{alignat}{2} \label{eq45.1a}
(\delta_t\epsilon_\phi^n,\epsilon_\phi^n)+\epsilon\kappa_1(\nabla\epsilon_\mu^n,\nabla\epsilon_\phi^n)&\leq({\Lambda}_h^n,\epsilon_\phi^n)
+(R^n_\phi,\epsilon_\phi^n),\\\label{eq45.1b}
-\kappa_1\epsilon(\nabla\epsilon_\phi^n,\nabla\epsilon_\mu^n)+\kappa_1(\epsilon_\mu^n,\epsilon_\mu^n)&=-\kappa_1(\widehat{\Lambda}_h^n,\epsilon_\mu^n).
\end{alignat}
\end{subequations}
Setting $\psi_h=\lambda\epsilon_\mu^n \in \cQ_h$,
$\theta_h=\lambda\delta_t\epsilon_\phi^n \in \cQ_h$, $\textbf{v}_h=\epsilon_\textbf{u}^n \in \cX_h$,
$\textbf{H}_h=\epsilon_\textbf{B}^n \in \cW_h$
in \eqref{eq37.1}, respectively,  it follows that
 \begin{subequations}\label{eq46.1}
\begin{alignat}{2} \label{eq46.1a}
\lambda(\delta_t\epsilon_\phi^n,\epsilon_\mu^n)+\lambda\epsilon\kappa_1\|\nabla\epsilon_\mu^n\|^2&\leq\lambda({\Lambda}_h^n,\epsilon_\mu^n)
+\lambda(R^n_\phi,\epsilon_\mu^n),\\\label{eq46.1b}
\epsilon\lambda(\nabla\epsilon_\phi^n,\nabla\delta_t\epsilon_\phi^n)-\lambda(\epsilon_\mu^n,\delta_t\epsilon_\phi^n)&=\lambda(\widehat{\Lambda}_h^n,\delta_t\epsilon_\phi^n),\\\label{eq46.1c}
 (\delta_t\epsilon_\textbf{u}^n,\epsilon_\textbf{u}^n)+c_0\nu_1\|\nabla \epsilon_\textbf{u}^n\|^2&\leq(\Phi_h^n,\epsilon_\textbf{u}^n)+(R_{\textbf{u}}^n,\epsilon_\textbf{u}^n),\\\label{eq46.1d}
S_c(\delta_t\epsilon_\textbf{B}^n,\epsilon_\textbf{B}^n)+c_1\eta_1S_c\|\nabla \epsilon_\textbf{B}^n\|^2&\leq(\widehat{\Phi}_h^n,\epsilon_\textbf{B}^n)+S_c(R_{\textbf{B}}^n,\epsilon_\textbf{B}^n).
\end{alignat}
\end{subequations}
Combining \eqref{eq45.1} and \eqref{eq46.1}, we obtain
 \begin{align}\label{eq47}
\frac{1}{2}\bigl(\|\epsilon_\phi^n\|^2&-\|\epsilon_\phi^{n-1}\|^2+\|\epsilon_\phi^n-\epsilon_\phi^{n-1}\|^2\bigr)\\&\nonumber\quad
+\frac{\epsilon\lambda}{2}\bigl(\|\nabla\epsilon_\phi^n\|^2-\|\nabla\epsilon_\phi^{n-1}\|^2+\|\nabla\epsilon_\phi^n-\nabla\epsilon_\phi^{n-1}\|^2\bigr)\\&\nonumber\quad
+\frac{1}{2}\bigl(\|\epsilon_\textbf{u}^n\|^2-\|\epsilon_\textbf{u}^{n-1}\|^2+\|\epsilon_\textbf{u}^n-\epsilon_\textbf{u}^{n-1}\|^2\bigr)\\&\nonumber\quad
+\frac{S_c}{2}\bigl(\|\epsilon_\textbf{B}^n\|^2-\|\epsilon_\textbf{B}^{n-1}\|^2+\|\epsilon_\textbf{B}^n-\epsilon_\textbf{B}^{n-1}\|^2\bigr)
+\Delta t\lambda\epsilon\kappa_1\|\nabla\epsilon_\mu^n\|^2
\\\nonumber&\quad
+\Delta t\kappa_1\|\epsilon_\mu^n\|^2+\Delta tc_0\nu_1\|\nabla \epsilon_\textbf{u}^n\|^2+\Delta tc_1\eta_1S_c\|\nabla \epsilon_\textbf{B}^n\|^2\\&\nonumber
\leq\Delta t\lambda(R^n_\phi,\epsilon_\mu^n)+\Delta t(R_\phi^n,\epsilon_\phi^n)+\Delta t(R_{\textbf{u}}^n,\epsilon_\textbf{u}^n)\\&\nonumber\quad
+\Delta tS_c(R_{\textbf{B}}^n,\epsilon_\textbf{B}^n)+\Delta t\lambda({\Lambda}_h^n,\epsilon_\mu^n)+\Delta t({\Lambda}_h^n,\epsilon_\phi^n)\\\nonumber&\quad
+\Delta t\lambda(\widehat{\Lambda}_h^n,\delta_t\epsilon_\phi^n)-\Delta t\kappa_1(\widehat{\Lambda}_h^n,\epsilon_\mu^n)\\\nonumber&\quad
+\Delta t(\Phi_h^n,\epsilon_\textbf{u}^n)
+\Delta t(\widehat{\Phi}_h^n,\epsilon_\textbf{B}^n)
\\&\nonumber
=\sum_{i=1}^{10}\Upsilon_i.
 \end{align}
We now bound the terms on the RHS of \eqref{eq47}.
  Define time dependent spatial mass average as follows
   \begin{align*}
 \overline{\epsilon}_\mu^n:=|\Omega|^{-1}( {\epsilon}_\mu^n,1), \ \qquad 1\leq n\leq N.
  \end{align*}
  For terms $\Upsilon_1-\Upsilon_4$, using \eqref{eq6.1}, Taylor's
 theorem, \textbf{Assumption 5} and Young's inequality, one finds that
 \begin{align}\label{eq48}
\Upsilon_1&=\Delta t\lambda(R^n_\phi,\epsilon_\mu^n-\overline{\epsilon}_\mu^n)\\\nonumber
 &\leq \Delta t\lambda\|\partial_t\phi^n-\delta_tr_h{\phi}^n\|\|\nabla\epsilon_\mu^n\|\\\nonumber
 &\leq C\Delta t\Bigl\{{\Delta t}\int_{t-\Delta t}^t\|\partial_{ss}^2\phi(s)\|^2_{L^2}ds\\\nonumber&\quad
 +\frac{h^{k+1}}{\sqrt{\Delta t}}\int^t_{t-\Delta t}\|\partial_s\phi(s)\|^2_{H^{k+1}}ds\Bigr\}^{\frac{1}{2}}\|\nabla\epsilon_\mu^n\|\\\nonumber&
 \leq C\Delta t\bigl({\Delta t}^2 +h^{2k+2}\bigr)+\frac{\kappa_1\lambda\Delta t}{8}\|\nabla\epsilon_\mu^n\|^2,\\\label{eq49}
\Upsilon_2&\leq C\Delta t\Bigl\{{\Delta t}\int_{t-\Delta t}^t\|\partial_{ss}^2\phi(s)\|^2_{L^2}ds\\\nonumber&\quad
 +\frac{h^{k+1}}{\sqrt{\Delta t}}\int^t_{t-\Delta t}\|\partial_s\phi(s)\|^2_{H^{k+1}}ds\Bigr\}^{\frac{1}{2}}\|\epsilon_\phi^n\|\\\nonumber&
  \leq C\Delta t\bigl({\Delta t}^2 +h^{2k+2}\bigr)+C\Delta t\|\epsilon_\phi^n\|^2.
 \end{align}
  Similarly we have
\begin{align}\label{eq50}
\Upsilon_3 \leq &C\Delta t\|\nabla \epsilon_\textbf{u}^n\|\|\partial_t\textbf{u}^n-\delta_tR_h{\textbf{u}}^n\| \\\nonumber
 \leq &C\Delta t\Bigl\{{\Delta t}\int_{t-\Delta t}^t\|\partial_{ss}^2\textbf{u}(s)\|^2_{L^2}ds\\\nonumber&\quad
 +\frac{h^{k+1}}{\sqrt{\Delta t}}\int^t_{t-\Delta t}\|(\partial_s\textbf{u}(s),\partial_sp(s))\|^2_{H^{k+1}\times H^{k})}ds\Bigr\}^{\frac{1}{2}}\|\nabla \epsilon_\textbf{u}^n\|^2\\\nonumber
 \leq& C\Delta t\bigl({\Delta t}^2 +{h^{2k+2}}\bigr)+\frac{c_0\nu_1\Delta t}{8}\|\nabla \epsilon_\textbf{u}^n\|^2,
 \end{align}
and
 \begin{align}\label{eq51}
\Upsilon_4 \leq &C\Delta t\|\nabla \epsilon_\textbf{B}^n\|\|\partial_t\textbf{B}^n-\delta_tR_{mh}{\textbf{B}}^n\| \\\nonumber
 \leq &C\Delta t\Bigl\{{\Delta t}\int_{t-\Delta t}^t\|\partial_{ss}^2\textbf{B}\|^2_{L^2}ds\\\nonumber&\quad
 +\frac{h^{k+1}}{\sqrt{\Delta t}}\int^t_{t-\Delta t}\|\partial_t\textbf{B}\|^2_{H^{k+1}}\Bigr\}^{\frac{1}{2}}\|\nabla \epsilon_\textbf{B}^n\|
 \\\nonumber
 \leq &C\Delta t\bigl({\Delta t}^2 +{h^{2k+2}}\bigr)+\frac{\Delta tc_1\eta_1S_c}{4}\|\nabla \epsilon_\textbf{B}^n\|^2.
 \end{align}
 For nonlinear term $\Upsilon_5$, adding and subtracting some terms, we can rewrite
\begin{align}\label{eq52}
\Upsilon_5&=\Delta t\lambda\Bigl\{ \epsilon a_\phi(\kappa^n(\phi^n)-\kappa^n(\phi^{n-1}); \mu^n, \epsilon_\mu^n)\\&\nonumber\quad
+ \epsilon a_\phi(\kappa^n(\phi^{n-1})-\kappa^n(\phi^{n-1}_h); \mu^n, \epsilon_\mu^n)+(\textbf{u}^{n}\cdot\nabla(\phi^{n}-\phi^{n-1}),\epsilon_\mu^n)\\&\nonumber
 \quad+\bigl((\textbf{u}^{n}-R_h\textbf{u}^{n})\cdot\nabla \phi^{n-1},\epsilon_\mu^n\bigr)
 +\bigl(R_h\textbf{u}^{n}\cdot\nabla (\phi^{n-1}-r_h\phi^{n-1}),\epsilon_\mu^n\bigr)\\&\nonumber
 \quad+(R_h\textbf{u}^{n}\cdot\nabla \epsilon_\phi^{n-1},\epsilon_\mu^n)
 +(\epsilon_\textbf{u}^{n}\cdot\nabla\phi_h^{n-1},\epsilon_\mu^n)\Bigr\}\\&\nonumber
 =\sum^{6}_{i=1}{\Lambda}_i^n +\Delta t\lambda(\epsilon_\textbf{u}^{n}\cdot\nabla\phi_h^{n-1},\epsilon_\mu^n).
 \end{align}
 For terms $\Lambda_1^n-\Lambda_2^n$, thanks to H\"{o}lder inequality,  Taylor's
 theorem and \eqref{eq33}, one finds that
  \begin{align*}
\Lambda_1^n+\Lambda_2^n&\leq C\Delta t|\kappa|_{C^{0,1}(\overline{\Omega}\times R;R)}\|\nabla\mu^n\|_{L^\infty}\Bigl\{\|\Delta t\delta_t\phi^n\|\\&\quad
+\|\phi^{n-1}-\widehat{r}_h\phi^{n-1}\|+\|\widehat{r}_h\phi^{n-1}-\phi^{n-1}_h\|\Bigr\}\|\nabla \epsilon_\mu^n\|\\&
\leq C\Delta t\Bigl\{\Big({\Delta t}\int^t_{t-\Delta t}\|\partial_s\phi(s)\|^2_{L^2}ds\Bigr)^{\frac{1}{2}}
+{h^{k+1}}\|\partial_t\phi\|_{H^{k+1}}\\&\quad
+\|\epsilon_\phi^{n-1}\|\Bigr\}\|\nabla \epsilon_\mu^n\|,
\end{align*}
where
  \begin{align*}
|\kappa|_{C^{0,1}(\overline{\Omega}\times R;R)}:=\sup\Bigl\{\frac{|\kappa(\textbf{x},\phi)
-\nu(\textbf{y},\varphi)|}{|(\textbf{x},\phi)-(\textbf{y},\varphi)|}: (\textbf{x},\phi),(\textbf{y},\varphi)\in \overline{\Omega}\times\R\Bigr\}.
  \end{align*}
  Using \eqref{eq6.1}, \eqref{eq29} and \eqref{eq33}, the terms $\Lambda_3^n-\Lambda_6^n$ can be bounded by
  \begin{align*}
  \Lambda_{3}^n &\leq C\Delta t\|\textbf{u}^{n}\|_{L^4}\|\Delta t\nabla\delta_t\phi^n\|\|\epsilon_\mu^n-\overline{\epsilon}_\mu^n\|_{L^4}\\&
\leq C\Delta t\Big({\Delta t}\int^t_{t-\Delta t}\|\partial_s\nabla\phi(s)\|^2_{L^2}ds\Bigr)^{\frac{1}{2}}\|\nabla\epsilon_\phi^n\|_{L^2},\\
\Lambda_{4}^n &\leq C\Delta t\|\textbf{u}^{n}-R_h\textbf{u}^{n}\|\nabla\phi^{n-1}\|_{L^4}\|\epsilon_\mu^n-\overline{\epsilon}_\mu^n\|_{L^4}\\&
\leq C\Delta t{h^{k+1}}\|(\textbf{u},p)\|_{H^{k+1}\times H^{k}}\|\nabla\epsilon_\mu^n\|,\\
\Lambda_{5}^n &\leq C\Delta t\|R_h\textbf{u}^{n}\|_{L^4}\|\nabla(\phi^{n-1}-r_h\phi^{n-1})\|\|\epsilon_\mu^n-\overline{\epsilon}_\mu^n\|_{L^4}\\&
\leq C\Delta t{h^{k+1}}\|\phi\|_{H^{k+2}}\|\nabla\epsilon_\phi^n\|,\\
\Lambda_{6}^n &\leq C\Delta t\|R_h\textbf{u}^{n}\|_{L^4}\|\nabla\epsilon_\phi^{n-1}\|\|\epsilon_\mu^n-\overline{\epsilon}_\mu^n\|_{L^4}
\leq C\Delta t\|\nabla\epsilon_\phi^{n-1}\|\|\nabla\epsilon_\mu^n\|.
\end{align*}
Combining \eqref{eq52} with above estimates, it follows that
  \begin{align}\label{eq53}
\Upsilon_5&\leq  C \Delta t\bigl(h^{2k+2}+\Delta t^2\bigr)+\Delta t\|\epsilon_\phi^{n-1}\|^2\\\nonumber&\quad
+\frac{\Delta t\lambda\kappa_1}{8}\|\nabla\epsilon_\mu^n\|^2+\Delta t\lambda(\epsilon_\textbf{u}^{n}\cdot\nabla\phi_h^{n-1},\epsilon_\mu^n).
\end{align}
Similarly, for nonlinear term $\Upsilon_6$, we rewrite
\begin{align}\label{eq54}
\Upsilon_6&=\Delta t\Bigl\{ \epsilon a_\phi(\kappa^n(\phi^n)-\kappa^n(\phi^{n-1}); \mu^n, \epsilon_\phi^n)\\&\nonumber\quad
+ \epsilon a_\phi(\kappa^n(\phi^{n-1})-\kappa^n(\phi^{n-1}_h); \mu^n, \epsilon_\phi^n)+(\textbf{u}^{n}\cdot\nabla(\phi^{n}-\phi^{n-1}),\epsilon_\phi^n)\\&\nonumber
 \quad+\bigl((\textbf{u}^{n}-R_h\textbf{u}^{n})\cdot\nabla \phi^{n-1},\epsilon_\phi^n\bigr)
 +\bigl(R_h\textbf{u}^{n}\cdot\nabla (\phi^{n-1}-r_h\phi^{n-1}),\epsilon_\phi^n\bigr)\\&\nonumber
 \quad+(R_h\textbf{u}^{n}\cdot\nabla \epsilon_\phi^{n-1},\epsilon_\phi^n)
 +(\epsilon_\textbf{u}^{n}\cdot\nabla\phi_h^{n-1},\epsilon_\phi^n)\Bigr\}\\&\nonumber
 =\sum^{13}_{i=7}{\Lambda}_i^n.
 \end{align}
  Applying \eqref{eq6.1}, \eqref{eq29} and \eqref{eq33}, we bound the terms $\Lambda_7^n-\Lambda_{13}^n$ as follows
  \begin{align*}
  \Lambda_7^n+\Lambda_8^n&\leq C{\Delta t}\Bigl\{\Big({\Delta t}\int^t_{t-\Delta t}\|\partial_s\phi(s)\|^2_{L^2}ds\Bigr)^{\frac{1}{2}}\\&\quad
+{h^{k+1}}
+\|\epsilon_\phi^{n-1}\|\Bigr\}\|\nabla\epsilon_\phi^{n}\|,\\
  \Lambda_{9}^n &\leq C\Delta t\|\textbf{u}^{n}\|_{L^\infty}\|\nabla(\phi^n-\phi^{n-1})\|\|\epsilon_\phi^n\|\\&
\leq C\Delta t\Big({\Delta t}\int^t_{t-\Delta t}\|\partial_s\nabla\phi(s)\|^2_{L^2}ds\Bigr)^{\frac{1}{2}}\|\epsilon_\phi^n\|,\\
\Lambda_{10}^n &\leq C\Delta t\|\textbf{u}^{n}-R_h\textbf{u}^{n}\|\nabla\phi^{n-1}\|_{L^\infty}\|\epsilon_\phi^n\|\\&
\leq C\Delta t{h^{k+1}}\|(\textbf{u},p)\|_{H^{k+1}\times H^{k}}\|\epsilon_\phi^n\|,\\
\Lambda_{11}^n &\leq C\Delta t\|R_h\textbf{u}^{n}\|_{L^\infty}\|\nabla(\phi^{n-1}-r_h\phi^{n-1})\|\|\epsilon_\phi^n\|\\&
\leq C\Delta t{h^{k+1}}\|\phi\|_{H^{k+2}}\|\epsilon_\phi^n\|,\\
\Lambda_{12}^n &\leq C\Delta t\|R_h\textbf{u}^{n}\|_{L^\infty}\|\nabla\epsilon_\phi^{n-1}\|\|\epsilon_\phi^n\|
\leq C\Delta t\|\nabla\epsilon_\phi^{n-1}\|\|\epsilon_\phi^n\|,\\
\Lambda_{13}^n &\leq C\Delta t\|\epsilon_\textbf{u}^{n}\|_{L^6}\|\nabla\phi_h^{n-1}\|_{L^3}\|\epsilon^n_\phi\|
\leq C\Delta t\|\nabla\epsilon_\textbf{u}^{n}\|\|\epsilon_\phi^n\|.
\end{align*}
Combining \eqref{eq54} with above estimates, one finds that
  \begin{align}\label{eq55}
 \Upsilon_6&\leq  C \Delta t\bigl(h^{2k+2}+\Delta t^2\bigr)+C\Delta t\|\epsilon_\phi^{n-1}\|^2+C\Delta t\|\nabla\epsilon_\phi^{n-1}\|^2\\\nonumber&\quad
+C\Delta t\|\nabla\epsilon_\phi^{n}\|^2
+\frac{\Delta tc_0\nu_1}{8}\|\nabla\epsilon_\textbf{u}^n\|^2+C\Delta t\|\epsilon_\phi^n\|^2.
\end{align}
 For nonlinear term $\Upsilon_7$, we can rewrite
\begin{align}\label{eq56}
\Upsilon_7&:=\Delta t\Bigl\{\lambda(\eta_\mu^n,\delta_t\epsilon_\phi^n)-\frac{\Delta t\lambda}{\epsilon}(\delta_t\phi^n,\delta_t\epsilon_\phi^n)
\\&\nonumber\quad-\lambda\epsilon^{-1}(e_\phi^{n-1},\delta_t\epsilon_\phi^n)+\lambda\epsilon^{-1}\bigl((\phi^{n})^3-(\phi^{n}_h)^3,\delta_t\epsilon_\phi^n\bigr)\Bigr\}\\&\nonumber
 =\Delta t\sum^{4}_{i=1}\widehat{{\Lambda}}_i^n.
 \end{align}
  Using \eqref{eq6.1}, \eqref{eq35} and \textbf{Lemmas 2-3}, we bound the terms $\widehat{\Lambda}_1^n-\widehat{\Lambda}_4^n$ as follows
  \begin{align*}
  \widehat{\Lambda}_{1}^n &\leq C\|\nabla\eta_\mu^n\|\|\delta_t\epsilon_\phi^n\|_{-1,h}\leq Ch^{k+1}\|\mu\|_{k+2}\|\delta_t\epsilon_\phi^n\|_{-1,h},\\
\widehat{\Lambda}_{2}^n &\leq \frac{1}{\epsilon}\|\Delta t\nabla\delta_t\phi^n\|\|\delta_t\epsilon_\phi^n\|_{-1,h}\leq C\Big({\Delta t}\int^t_{t-\Delta t}\|\partial_s\nabla\phi(s)\|^2_{L^2}ds\Bigr)^{\frac{1}{2}}\|\delta_t\epsilon_\phi^n\|_{-1,h},\\
\widehat{\Lambda}_{3}^n &\leq C\|e_\phi^{n-1}\|\|\delta_t\epsilon_\phi^n\|_{-1,h}
\leq C\|\epsilon_\phi^{n-1}\|\|\delta_t\epsilon_\phi^n\|_{-1,h}+C\|\eta_\phi^{n-1}\|\|\delta_t\epsilon_\phi^n\|_{-1,h},\\
\widehat{\Lambda}_{4}^n &\leq C\|\nabla(\phi^{n})^3-(\phi^{n}_h)^3)\|\|\delta_t\epsilon_\phi^n\|_{-1,h}\leq C\|\nabla(\phi^{n}-\phi^{n}_h)\|\|\delta_t\epsilon_\phi^n\|_{-1,h}\\&
\leq C\|\nabla\epsilon_\phi^{n}\|\|\delta_t\epsilon_\phi^n\|_{-1,h}+C\|\nabla\eta_\phi^{n}\|\|\delta_t\epsilon_\phi^n\|_{-1,h}.
\end{align*}
Combining \eqref{eq56} with above estimates, it follows that
  \begin{align}\label{eq57}
\Upsilon_7&\leq  C\Delta t\bigl(h^{2k+2}+\Delta t^2\bigr)
+C\Delta t\|\epsilon_\phi^{n-1}\|^2\\&\nonumber\quad
+C\Delta t\|\nabla\epsilon_\phi^{n}\|^2+\alpha\Delta t\|\delta_t\epsilon_\phi^n\|^2_{-1,h}.
\end{align}
Taking $T_h(\delta_t\epsilon_\phi^n)$ in \eqref{eq37.1a}, and using \eqref{eq6.1}, \textbf{Lemmas 1-2}, we derive
 \begin{align}\label{eq58}
\|\delta_t\epsilon_\phi^n\|^2_{-1,h}&=-\epsilon a_\phi(\kappa^n(\phi^{n-1}_h);\epsilon_\mu^n,T_h(\delta_t\epsilon_\phi^n))+(R^n_\phi,T_h(\delta_t\epsilon_\phi^n))+({\Lambda}_h^n,T_h(\delta_t\epsilon_\phi^n))\\&\nonumber
\leq \epsilon\kappa_2\|\nabla\epsilon_\mu^n\|\|\delta_t\epsilon_\phi^n\|_{-1,h}+\|\partial_t\phi^n-\delta_tr_h{\phi}^n\|\|T_h(\delta_t\epsilon_\phi^n)\|\\&\nonumber\quad+
\|\textbf{u}^{n}\|_{L^3}\|\Delta t\nabla\delta_t\phi^{n}\|\|T_h(\delta_t\epsilon_\phi^n)\|_{L^6}
+\|\eta_\textbf{u}^{n}\|_{L^3}\|\nabla \phi^{n-1}\|\|T_h(\delta_t\epsilon_\phi^n)\|_{L^6}\\&\nonumber\quad
 +\|R_h\textbf{u}^{n}\|_{L^3}\|\nabla \eta_\phi^{n-1}\|\|T_h(\delta_t\epsilon_\phi^n)\|_{L^6}
 +\|R_h\textbf{u}^{n}\|_{L^3}\|\nabla \epsilon_\phi^{n-1}\|\|T_h(\delta_t\epsilon_\phi^n)\|_{L^6}\\&\nonumber\quad
 +\|\epsilon_\textbf{u}^{n}\|_{L^3}\|\nabla\phi_h^{n-1}\|\|T_h(\delta_t\epsilon_\phi^n)\|_{L^6}\\&\nonumber
 \leq\frac{1}{2}\|\delta_t\epsilon_\phi^n\|^2_{-1,h}+{C_1}\|\nabla\epsilon_\mu^n\|^2
+C_2\|\nabla\epsilon_\textbf{u}^{n}\|^2\\&\nonumber\quad
+C\|\nabla\epsilon_\phi^{n-1}\|^2+C\bigl(h^{2k+2}+\Delta t^2\bigr).
 \end{align}
Then we obtain
  \begin{align}\label{eq59}
\|\delta_t\epsilon_\phi^n\|^2_{-1,h}
 \leq 2{C_1}\|\nabla\epsilon_\mu^n\|^2
+2C_2\|\nabla\epsilon_\textbf{u}^{n}\|^2+C\|\nabla\epsilon_\phi^{n-1}\|^2+C\bigl(h^{2k+2}+\Delta t^2\bigr).
 \end{align}
  From \eqref{eq57} and \eqref{eq59}, and taking $\alpha=\min\bigl\{\frac{\lambda\kappa_1}{8C_1},\frac{c_0\nu_1}{8C_2}\bigr\}$ in \eqref{eq57}, one finds that
  \begin{align}\label{eq60}
\Upsilon_7&\leq  C\Delta t\bigl(h^{2k+2}+\Delta t^2\bigr)
+\frac{\Delta tc_0\nu_1}{8}\|\nabla\epsilon_\textbf{u}^{n}\|^2+\frac{\Delta t\lambda\kappa_1}{8}\|\nabla\epsilon_\mu^n\|^2\\&\nonumber\quad
+C\Delta t\|\nabla\epsilon_\phi^{n-1}\|^2+C\Delta t\|\nabla\epsilon_\phi^{n}\|^2.
 \end{align}
 Similarly, for nonlinear term $\Upsilon_8$, we rewrite
\begin{align}\label{eq61}
\Upsilon_8&:=\Delta t\Bigl\{(\eta_\mu^n,\epsilon_\mu^n)+\frac{\Delta t}{\epsilon}(\delta_t\phi^n,\epsilon_\mu^n)
\\&\nonumber\quad+\epsilon^{-1}(e_\phi^{n-1},\epsilon_\mu^n)+\epsilon^{-1}\bigl((\phi^{n})^3-(\phi^{n}_h)^3,\epsilon_\mu^n\bigr)\Bigr\}\\&\nonumber
 =\sum^{8}_{i=5}\widehat{{\Lambda}}_i^n.
 \end{align}
  Using \eqref{eq6.1}, \eqref{eq33}, \eqref{eq35} and \textbf{Lemma 3}, we bound the terms $\widehat{\Lambda}_5^n-\widehat{\Lambda}_8^n$ as follows
  \begin{align*}
  \widehat{\Lambda}_{5}^n &\leq C\Delta t\|\eta_\mu^n\|\|\epsilon_\mu^n\|\leq Ch^{k+1}\|\mu\|_{k+1}\|\epsilon_\mu^n\|,\\
\widehat{\Lambda}_{6}^n &\leq \Delta t\frac{1}{\epsilon}\|\Delta t\nabla\delta_t\phi^n\|\|\epsilon_\mu^n\|
\leq C\Delta t\Big({\Delta t}\int^t_{t-\Delta t}\|\partial_s\nabla\phi(s)\|^2_{L^2}ds\Bigr)^{\frac{1}{2}}\|\epsilon_\mu^n\|,\\
\widehat{\Lambda}_{7}^n &\leq C\Delta t\|e_\phi^{n-1}\|\|\epsilon_\mu^n\|\leq C\Delta t\|\epsilon_\phi^{n-1}\|\|\epsilon_\mu^n\|+C\Delta t\|\eta_\phi^{n-1}\|\|\epsilon_\mu^n\|\\
\widehat{\Lambda}_{8}^n &\leq C\Delta t\|(\phi^{n})^3-(\phi^{n}_h)^3)\|\|\epsilon_\mu^n\|\leq C\Delta t\|(\phi^{n}-\phi^{n}_h)\|\|\epsilon_\mu^n\|\\&
\leq C\Delta t\|\epsilon_\phi^{n}\|\|\epsilon_\mu^n\|+C\Delta t\|\nabla\eta_\phi^{n}\|\|\epsilon_\mu^n\|.
\end{align*}
Combining \eqref{eq61} with above inequalities, we have
  \begin{align}\label{eq62}
\Upsilon_8&\leq  C\Delta t\bigl(h^{2k+2}+\Delta t^2\bigr)+\frac{\Delta t\kappa_1}{8}\|\epsilon_\mu^n\|^2\\&\nonumber\quad
+C\Delta t\|\epsilon_\phi^{n}\|^2+C\Delta t\|\epsilon_\phi^{n-1}\|^2.
\end{align}
For nonlinear term $\Upsilon_9$, by adding and
subtracting some terms, we  rewrite
 \begin{align}\label{eq63}
\Upsilon_9&=\Delta t\Big\{ a_f\bigl(\nu^n(\phi^n)-\nu^n(\phi^{n-1}); \textbf{u}^n, \epsilon_\textbf{u}^n\bigr)\\&\nonumber\quad
 +a_f\bigl(\nu^n(\phi^{n-1})-\nu^n(\phi^{n-1}_h); \textbf{u}^n, \epsilon_\textbf{u}^n\bigr)+b(\eta_\textbf{u}^n,\textbf{u}^n,\epsilon_\textbf{u}^n)\\&\nonumber
 \quad+b(R_h\textbf{u}^n-R_h\textbf{u}^{n-1},\textbf{u}^n,\epsilon_\textbf{u}^n)+b(\epsilon_\textbf{u}^{n-1},\textbf{u}^n,\epsilon_\textbf{u}^n)\\&\nonumber
 \quad+b(\textbf{u}_h^{n-1},\eta_\textbf{u}^n,\epsilon_\textbf{u}^n)+b(\textbf{u}_h^{n-1},\epsilon_\textbf{u}^n,\epsilon_\textbf{u}^n)+S_cc_{\widehat{B}}(\eta_\textbf{B}^n,\textbf{B}^n,\epsilon_\textbf{u}^n)\\&\nonumber
  \quad+S_cc_{\widehat{B}}(R_{mh}\textbf{B}^n-R_{mh}\textbf{B}^{n-1},\textbf{B}^n,\epsilon_\textbf{u}^n)+S_cc_{\widehat{B}}(\epsilon_\textbf{B}^{n-1},\textbf{B}^n,\epsilon_\textbf{u}^n)\\&\nonumber
   \quad+S_cc_{\widehat{B}}(\textbf{B}_h^{n-1},\eta_\textbf{B}^n,\epsilon_\textbf{u}^n)+S_cc_{\widehat{B}}(\textbf{B}_h^{n-1},\epsilon_\textbf{B}^n,\epsilon_\textbf{u}^n)\\&\nonumber\quad
   -\lambda(\mu^n\nabla\eta_\phi^n,\epsilon_\textbf{u}^n)-\lambda\bigl(\mu^n\nabla(r_h\phi^n-r_h\phi^{n-1}),\epsilon_\textbf{u}^n\bigr)\\&\nonumber\quad
    -\lambda(\mu^n\nabla\epsilon_\phi^{n-1},\epsilon_\textbf{u}^n)-\lambda(\eta_\mu^n\nabla\phi_h^{n-1},\epsilon_\textbf{u}^n)-\lambda(\epsilon_\mu^n\nabla\phi_h^{n-1},\epsilon_\textbf{u}^n)\\&\nonumber
  =\Delta t\Big\{\sum^{14}_{i=1}\Phi_i^n +b(\textbf{u}_h^{n-1},\epsilon_\textbf{u}^n,\epsilon_\textbf{u}^n)+S_cc_{\widehat{B}}(\textbf{B}_h^{n-1},\epsilon_\textbf{B}^n,\epsilon_\textbf{u}^n)
  \\&\nonumber\quad -\lambda(\epsilon_\mu^n\nabla\phi_h^{n-1},\epsilon_\textbf{u}^n) \Bigr\}.
 \end{align}
 For the terms $\Phi_1^n-\Phi_2^n$, using H\"{o}lder inequality, \eqref{eq35} and Taylor's
 theorem, we have
  \begin{align*}
\Phi_1^n+\Phi_2^n&\leq C|\nu|_{C^{0,1}(\overline{\Omega}\times R;R)}\|\mathbb{D}(\textbf{u}^n)\|_{L^\infty}\Bigl\{\|\phi^n-\phi^{n-1}\|\\&\quad
+\|\phi^{n-1}-\widehat{r}_h\phi^{n-1}\|+\|\widehat{r}_h\phi^{n-1}-\phi^{n-1}_h\|\Bigr\}\|\nabla \epsilon_\textbf{u}^n\|\\&
\leq C\Bigl\{\Big({\Delta t}\int^t_{t-\Delta t}\|\partial_s\phi(s)\|^2_{L^2}ds\Bigr)^{\frac{1}{2}}
+{h^{k+1}}\\&\quad
+\|\epsilon_\phi^{n-1}\|\Bigr\}\|\nabla \epsilon_\textbf{u}^n\|,
\end{align*}
Making use of \eqref{eq6.1} and \eqref{eq29}, we estimate the terms $\Phi_3^n-\Phi_6^n$  as follows
  \begin{align*}
\Phi_{3}^n &\leq C\|\nabla\eta_\textbf{u}^n\|\|\nabla \textbf{u}^n\|\|\nabla \epsilon_\textbf{u}^n\|\\&
\leq C{h^{k+1}}\|(\textbf{u},p)\|_{H^{k+2}\times H^{k+1}}\|\nabla \epsilon_\textbf{u}^n\|,\\
\Phi_{4}^n &\leq C\|\Delta tR_{mh}\delta_t\textbf{B}^n\|\|\nabla \textbf{B}^n\|_{L^3}\|\nabla \epsilon_\textbf{u}^n\|\\&
\leq C\Big({\Delta t}\int^t_{t-\Delta t}\|\partial_s\textbf{B}(s)\|^2_{L^2}ds\Bigr)^{\frac{1}{2}}\|\nabla \epsilon_\textbf{u}^n\|,\\
\Phi_{5}^n &\leq C\|\epsilon_\textbf{u}^{n-1}\|\|\nabla \textbf{u}^n\|_{L^3}\|\nabla \epsilon_\textbf{u}^n\|\\&
\leq\|\epsilon_\textbf{u}^{n-1}\|\|\nabla \epsilon_\textbf{u}^n\|,\\
\Phi_{6}^n &\leq C\|\nabla\textbf{u}_h^{n-1}\|\|\nabla\eta_\textbf{u}^n\|\|\nabla\epsilon_\textbf{u}^n\|\\&
\leq C{h^{k+1}}\|(\textbf{u},p)\|_{H^{k+2}\times H^{k+1}}\|\nabla \epsilon_\textbf{u}^n\|.
\end{align*}
Thanks to \eqref{eq6.1} and \eqref{eq31},  the terms $\Phi_7^n-\Phi_{10}^n$ can be bounded by
  \begin{align*}
\Phi_{7}^n &\leq C\|\eta_\textbf{B}^n\|\|\textbf{B}^n\|_{L^3}\|\nabla \epsilon_\textbf{u}^n\|\\&
\leq C{h^{k+1}}\|\textbf{B}\|_{H^{k+1}}\|\nabla \epsilon_\textbf{u}^n\|,\\
\Phi_{8}^n &\leq C\|\Delta tR_{mh}\delta_t\textbf{B}^n\|\|\curl\textbf{B}^n\|_{L^3}\|\nabla \epsilon_\textbf{u}^n\|\\&
\leq C\Big({\Delta t}\int^t_{t-\Delta t}\|\partial_s\textbf{B}(s)\|^2_{L^2}ds\Bigr)^{\frac{1}{2}}\|\nabla \epsilon_\textbf{u}^n\|,\\
\Phi_{9}^n &\leq C\|\epsilon_\textbf{B}^{n-1}\|\|\curl\textbf{B}^n\|_{L^3}\|\nabla \epsilon_\textbf{u}^n\|
\leq C\|\epsilon_\textbf{B}^{n-1}\|\|\nabla \epsilon_\textbf{u}^n\|,\\
\Phi_{10}^n &\leq C\|\eta_\textbf{B}^{n-1}\|\|\curl\textbf{B}^n\|_{L^3}\|\nabla \epsilon_\textbf{u}^n\|\\&
\leq C{h^{k+1}}\|\textbf{B}\|_{H^{k+1}}\|\nabla \epsilon_\textbf{u}^n\|.
\end{align*}
Applying H\"{o}lder inequality, \eqref{eq6.1}, \eqref{eq33} and \eqref{eq35}, the terms $\Phi_{11}^n-\Phi_{14}^n$ can be estimated as follows
  \begin{align*}
\Phi_{11}^n &\leq C\|\mu^n\|_{L^3}\|\nabla\eta_\phi^n\|\|\epsilon_\textbf{u}^n\|_{L^6}\\&
\leq C{h^{k+1}}\|\phi\|_{H^{k+2}}\|\nabla \epsilon_\textbf{u}^n\|,\\
\Phi_{12}^n &\leq C\|\mu^n\|_{L^3}\|\Delta t\nabla \delta_tr_h\phi^n\|\|\epsilon_\textbf{u}^n\|_{L^6}\\&
\leq C\Big({\Delta t}\int^t_{t-\Delta t}\|\partial_s\nabla\phi(s)\|^2_{L^2}ds\Bigr)^{\frac{1}{2}}\|\nabla \epsilon_\textbf{u}^n\|,\\
\Phi_{13}^n &\leq C\|\mu^n\|_{L^3}\|\nabla\epsilon_\phi^{n-1}\|\|\epsilon_\textbf{u}^n\|_{L^6}\\&
\leq C\|\nabla\epsilon_\phi^{n-1}\|\|\nabla \epsilon_\textbf{u}^n\|,\\
\Phi_{14}^n &\leq C\|\eta_\mu^n\|\|\nabla\phi_h^{n-1}\|\|\epsilon_\textbf{u}^n\|_{L^6}\\&
\leq C{h^{k+1}}\|\mu\|_{H^{k+1}}\|\nabla \epsilon_\textbf{u}^n\|.
\end{align*}
Combining \eqref{eq63} with above bounds, one finds that
  \begin{align}\label{eq64}
\Upsilon_9&\leq  C\Delta t\Bigl\{h^{k+1}+\Delta t+\|\epsilon_{\textbf{u}}^{n-1}\|+\|\epsilon_{\textbf{B}}^{n-1}\|+\|\epsilon_\phi^{n-1}\|\Bigr\}\|\nabla \epsilon_\textbf{u}^n\|\\&\nonumber\quad
+\Delta tb(\textbf{u}_h^{n-1},\epsilon_\textbf{u}^n,\epsilon_\textbf{u}^n) +\Delta tS_cc_{\widehat{B}}(\textbf{B}_h^{n-1},\epsilon_\textbf{B}^n,\epsilon_\textbf{u}^n)-\Delta t\lambda(\epsilon_\mu^n\nabla\phi_h^{n-1},\epsilon_\textbf{u}^n)\\\nonumber&
  \leq  C\Bigl\{\Delta t\bigl(h^{2k+2}+\Delta t^2\bigr)+\Delta t\|\epsilon_{\textbf{u}}^{n-1}\|^2+\Delta t\|\epsilon_{\textbf{B}}^{n-1}\|^2+\Delta t\|\epsilon_\phi^{n-1}\|^2\Bigr\}+\frac{\Delta tc_0 \nu_1}{8}\|\nabla \epsilon_\textbf{u}^n\|^2\\\nonumber&
  \quad+\Delta tb(\textbf{u}_h^{n-1},\epsilon_\textbf{u}^n,\epsilon_\textbf{u}^n)+\Delta tS_cc_{\widehat{B}}(\textbf{B}_h^{n-1},\epsilon_\textbf{B}^n,\epsilon_\textbf{u}^n)-\Delta t\lambda(\epsilon_\mu^n\nabla\phi_h^{n-1},\epsilon_\textbf{u}^n).
\end{align}
For nonlinear term $\Upsilon_{10}$, by adding and
subtracting some terms, we can rewrite
 \begin{align}\label{eq65}
 \Upsilon_{10}&=\Delta t\Bigl\{S_ca_\textbf{B}(\eta^n(\phi^n)-\eta^n(\phi^{n-1}); \textbf{B}^n, \epsilon_\textbf{B}^n)\\&\nonumber\quad
+ S_ca_\textbf{B}(\eta^n(\phi^{n-1})-\eta^n(\phi^{n-1}_h); \textbf{B}^n, \epsilon_\textbf{B}^n)
-S_cc_{\widetilde{B}}(\textbf{u}^{n},\eta_\textbf{B}^{n},\epsilon_\textbf{B}^n)\\&\nonumber
 \quad-S_cc_{\widetilde{B}}(\textbf{u}^{n},R_{mh}\textbf{B}^{n}-R_{mh}\textbf{B}^{n-1},\epsilon_\textbf{B}^n)
 -S_cc_{\widetilde{B}}(\textbf{u}^{n},\epsilon_\textbf{B}^{n-1},\epsilon_\textbf{B}^n)\\&\nonumber
 \quad-S_cc_{\widetilde{B}}(\eta_\textbf{u}^{n},\textbf{B}_h^{n-1},\epsilon_\textbf{B}^n)
-S_cc_{\widetilde{B}}(\epsilon_\textbf{u}^{n},\textbf{B}_h^{n-1},\epsilon_\textbf{B}^n)\Bigr\}\\&\nonumber
 =\Delta t\Bigl\{\sum^{6}_{i=1}\widehat{\Phi}_i^n-S_cc_{\widetilde{B}}(\epsilon_\textbf{u}^{n},\textbf{B}_h^{n-1},\epsilon_\textbf{B}^n)\Bigl\}.
 \end{align}
  For the terms $\widehat{\Phi}_1^n-\widehat{\Phi}_2^n$, thanks to H\"{o}lder inequality, \eqref{eq35} and Taylor's
 theorem, we get
   \begin{align*}
\widehat{\Phi}_1^n+\widehat{\Phi}_2^n&\leq C|\eta|_{C^{0,1}(\overline{\Omega}\times R;R)}\|\nabla\textbf{B}^n\|_{L^\infty}\Bigl\{\|\Delta t\delta_t\phi^n\|\\&\quad
+\|\phi^{n-1}-\widehat{r}_h\phi^{n-1}\|+\|\widehat{r}_h\phi^{n-1}-\phi^{n-1}_h\|\Bigr\}\|\nabla \epsilon_\textbf{B}^n\|\\&
\leq C\Bigl\{\Big({\Delta t}\int^t_{t-\Delta t}\|\partial_s\phi(s)\|^2_{L^2}ds\Bigr)^{\frac{1}{2}}
+{h^{k+1}}\\&\quad
+\|\epsilon_\phi^{n-1}\|\Bigr\}\|\nabla \epsilon_\textbf{B}^n\|
\end{align*}
By \eqref{eq6.1}, \eqref{eq29} and \eqref{eq31}, the terms $\widehat{\Phi}_3^n-\widehat{\Phi}_6^n$ can be bounded by
  \begin{align*}
\widehat{\Phi}_{3}^n &\leq C\|\textbf{u}^{n}\|_{L^{\infty}}\|\eta_\textbf{B}^{n}\|\|\curl\epsilon_\textbf{B}^n\|\\&
\leq C{h^{k+1}}\|\textbf{B}\|_{H^{k+1}}\|\nabla \epsilon_\textbf{B}^n\|,\\
\widehat{\Phi}_{4}^n &\leq C\|\textbf{u}^n\|_{L^{\infty}}\|\Delta tR_{mh}\delta_t\textbf{B}^{n}\|\|\curl\epsilon_\textbf{B}^n\|\\&
\leq C\Big({\Delta t}\int^t_{t-\Delta t}\|\partial_s\textbf{B}(s)\|^2_{L^2}ds\Bigr)^{\frac{1}{2}}\|\nabla \epsilon_\textbf{B}^n\|,\\
\widehat{\Phi}_{5}^n &\leq C\|\textbf{u}^n\|_{L^{\infty}}\|\epsilon_\textbf{B}^{n-1}\|\|\curl \epsilon_\textbf{B}^n\|\\&
\leq C\|\epsilon_\textbf{B}^{n-1}\|\|\nabla \epsilon_\textbf{B}^n\|,\\
\widehat{\Phi}_{6}^n &\leq C\|\eta_\textbf{u}^{n}\|_{L^{6}}\|\textbf{B}_h^{n-1}\|_{L^3}\|\curl\epsilon_\textbf{B}^n\|\\&
\leq C{h^{k+1}}\|(\textbf{u},p)\|_{H^{k+2}\times H^{k+1})}\|\textbf{B}_h^{n-1}\|_{L^3}\|\nabla \epsilon_\textbf{B}^n\|.
\end{align*}
Combining \eqref{eq65} with above inequalities, it follows that
   \begin{align}\label{eq66}
\Upsilon_{10}&\leq  C \Delta t\Bigl\{h^{k+1}+\Delta t+\|\epsilon_{\textbf{B}}^{n-1}\|+\|\epsilon_\phi^{n-1}\|\Bigr\}\|\nabla \epsilon_\textbf{B}^n\|\\\nonumber&
  \quad- \Delta tS_cc_{\widetilde{B}}(\epsilon_\textbf{u}^{n},\textbf{B}_h^{n-1},\epsilon_\textbf{B}^n)\\\nonumber&
 \leq  C\Delta t\bigl(h^{2k+2}+\Delta t^2\bigr)+\Delta t\|\epsilon_{\textbf{B}}^{n-1}\|^2+\Delta t\|\epsilon_\phi^{n-1}\|^2\\\nonumber&\quad+\frac{c_1\eta_1S_c\Delta t}{4}\|\nabla \epsilon_\textbf{B}^n\|^2
 - \Delta tS_cc_{\widetilde{B}}(\epsilon_\textbf{u}^{n},\textbf{B}_h^{n-1},\epsilon_\textbf{B}^n).
\end{align}
 Combining \eqref{eq47} with \eqref{eq7}-\eqref{eq8}, \eqref{eq48}-\eqref{eq51}, \eqref{eq53}, \eqref{eq55}, \eqref{eq60}, \eqref{eq62}, \eqref{eq64} and \eqref{eq66}, we have
  \begin{align}\label{eq67}
\frac{1}{2}\bigl(\|\epsilon_\phi^n\|^2&-\|\epsilon_\phi^{n-1}\|^2+\|\epsilon_\phi^n-\epsilon_\phi^{n-1}\|^2\bigr)\\&\nonumber\quad
+\frac{\epsilon\lambda}{2}\bigl(\|\nabla\epsilon_\phi^n\|^2-\|\nabla\epsilon_\phi^{n-1}\|^2+\|\nabla\epsilon_\phi^n-\nabla\epsilon_\phi^{n-1}\|^2\bigr)\\&\nonumber\quad
+\frac{1}{2}\bigl(\|\epsilon_\textbf{u}^n\|^2-\|\epsilon_\textbf{u}^{n-1}\|^2+\|\epsilon_\textbf{u}^n-\epsilon_\textbf{u}^{n-1}\|^2\bigr)\\&\nonumber\quad
+\frac{S_c}{2}\bigl(\|\epsilon_\textbf{B}^n\|^2-\|\epsilon_\textbf{B}^{n-1}\|^2+\|\epsilon_\textbf{B}^n-\epsilon_\textbf{B}^{n-1}\|^2\bigr)
+\Delta t\lambda\kappa_1\|\nabla\epsilon_\mu^n\|^2
\\\nonumber&\quad
+\Delta t\kappa_1\|\epsilon_\mu^n\|^2+\Delta tc_0\nu_1\|\nabla \epsilon_\textbf{u}^n\|^2+\Delta tc_1\eta_1S_c\|\nabla \epsilon_\textbf{B}^n\|^2\\&\nonumber
\leq C\Delta t\bigl(h^{2k+2}+\Delta t^2\bigr)+C\Delta t\bigl(\|\epsilon_{\phi}^{n}\|^2+\|\nabla\epsilon_\phi^{n}\|^2\bigr)\\&\nonumber\quad
+C\Delta t\bigl(\|\epsilon_{\textbf{u}}^{n-1}\|^2+\|\epsilon_{\textbf{B}}^{n-1}\|^2+\|\epsilon_\phi^{n-1}\|^2+\|\nabla\epsilon_\phi^{n-1}\|^2\bigr).
\end{align}
Summing \eqref{eq67} from $n=1$ to $m$ and using \eqref{eq36}, we obtain
 \begin{align}\label{eq68}
\|\epsilon_\phi^m\|^2&+{\epsilon}\lambda\|\nabla\epsilon_\phi^m\|^2+\|\epsilon_\textbf{u}^m\|^2+\|\epsilon_\textbf{B}^m\|^2
+\Delta t\lambda\kappa_1\sum_{n=1}^m\|\nabla\epsilon_\mu^n\|^2
\\\nonumber&\quad
+\Delta t\kappa_1\sum_{n=1}^m\|\epsilon_\mu^n\|^2+\Delta tc_0\nu_1\sum_{n=1}^m\|\nabla \epsilon_\textbf{u}^n\|^2+\Delta tc_1\eta_1S_c\sum_{n=1}^m\|\nabla \epsilon_\textbf{B}^n\|^2\\&\nonumber
\leq C\bigl(h^{2k+2}+\Delta t^2\bigr)+C_3\Delta t\sum_{n=1}^m\|\epsilon_{\phi}^{n}\|^2+C_4\Delta t\sum_{n=1}^m\|\nabla\epsilon_\phi^{n}\|^2\\&\nonumber\quad
+C_5\Delta t\sum_{n=1}^m\|\epsilon_{\textbf{u}}^{n-1}\|^2+C_6\Delta t\sum_{n=1}^m\|\epsilon_{\textbf{B}}^{n-1}\|^2\\&\nonumber\quad
+C_7\Delta t\sum_{n=1}^m\|\epsilon_\phi^{n-1}\|^2+C_8\Delta t\sum_{n=1}^m\|\nabla\epsilon_\phi^{n-1}\|^2.
\end{align}
If $0<\Delta t\leq \Delta t_0:=\frac{1}{2\max\{C_3,C_4\epsilon^{-1}\}}<\frac{1}{\max\{C_3,C_4\epsilon^{-1}\}}$, since $1\leq \frac{1}{1-\max\{C_3,C_4\epsilon^{-1}\}\Delta t}\leq 2$, it follows from \eqref{eq68} that
 \begin{align}\label{eq69}
\|\epsilon_\phi^m\|^2&+{\epsilon}\lambda\|\nabla\epsilon_\phi^m\|^2+\|\epsilon_\textbf{u}^m\|^2+\|\epsilon_\textbf{B}^m\|^2
+\Delta t\lambda\kappa_1\sum_{n=1}^m\|\nabla\epsilon_\mu^n\|^2
\\\nonumber&\quad
+\Delta t\kappa_1\sum_{n=1}^m\|\epsilon_\mu^n\|^2+\Delta tc_0\nu_1\sum_{n=1}^m\|\nabla \epsilon_\textbf{u}^n\|^2+\Delta tc_1\eta_1S_c\sum_{n=1}^m\|\nabla \epsilon_\textbf{B}^n\|^2\\&\nonumber
\leq C\bigl(h^{2k+2}+\Delta t^2\bigr)+\frac{(C_3+C_7)\Delta t}{1-\max\{C_3,C_4\epsilon^{-1}\}}\sum_{n=1}^{m}\|\epsilon_{\phi}^{n-1}\|^2\\&\nonumber\quad
+\frac{(C_4+C_8)\Delta t}{1-\max\{C_3,C_4\epsilon^{-1}\}}\sum_{n=1}^{m}\|\nabla\epsilon_\phi^{n-1}\|^2
+\frac{C_5\Delta t}{1-\max\{C_3,C_4\epsilon^{-1}\}}\sum_{n=1}^m\|\epsilon_{\textbf{u}}^{n-1}\|^2
\\&\nonumber\quad+\frac{C_6\Delta t}{1-\max\{C_3,C_4\epsilon^{-1}\}}\sum_{n=1}^m\|\epsilon_{\textbf{B}}^{n-1}\|^2.
 \end{align}
By using the discrete Gr\"{o}nwall inequality, one finds that
  \begin{align}\label{eq70}
\|\epsilon_\phi^m\|^2&+{\epsilon}\lambda\|\nabla\epsilon_\phi^m\|^2+\|\epsilon_\textbf{u}^m\|^2+\|\epsilon_\textbf{B}^m\|^2
+\Delta t\lambda\kappa_1\sum_{n=1}^m\|\nabla\epsilon_\mu^n\|^2
\\\nonumber&
+\Delta t\kappa_1\sum_{n=1}^m\|\epsilon_\mu^n\|^2+\Delta tc_0\nu_1\sum_{n=1}^m\|\nabla \epsilon_\textbf{u}^n\|^2+\Delta tc_1\eta_1S_c\sum_{n=1}^m\|\nabla \epsilon_\textbf{B}^n\|^2\\&\nonumber
\leq C\bigl(h^{2k+2}+\Delta t^2\bigr).
\end{align}
The desired result follow from \eqref{eq29}, \eqref{eq31}, \eqref{eq33}, \eqref{eq35} and the triangle inequality. The proof is completed.
 $$\eqno\Box$$

In order to derive optimal error estimates on pressure, we define the seminorm for $\textbf{v} \in H^{-1}(\Omega)^d$:
\begin{align*}
\|\textbf{v}\|_{*,h}=\sup_{\textbf{v}_h\in \cX_h}\frac{(\textbf{v},\textbf{v}_h)}{\|\nabla\textbf{v}_h\|}.
 \end{align*}
And we define the following inverse Stokes operator: $S:H^{-1}(\Omega)^d\rightarrow\cX$ and discrete inverse Stokes operator: $S_h:H^{-1}(\Omega)^d\rightarrow\cX_h$.
 Let $\hat{\nu} \in C^{0,1}(\overline{\Omega};R^{+})$ satisfy \textbf{Assumption A4},
  for all $\textbf{v} \in H^{-1}(\Omega)^d$, $(S(\textbf{v}),r)\in \cX\times\cM$ is the weak solution of the following problem:
 \begin{align*}
a_f(\widehat{\nu};S(\textbf{v}),\textbf{v})+d(\textbf{v},r)&=(\textbf{v},\textbf{v}),\ &&\forall\, \textbf{v}\in\cX,\\
d(S(\textbf{v}),w)&=0,\ &&\forall\, w\in\cM.
 \end{align*}
We have $H^2$ regularity results \cite{Temam1983}:
 \begin{align*}
\|S(\textbf{v})\|_2+\|\nabla r\|\leq C\|\textbf{v}\|, \forall\, \textbf{v}\in L^2(\Omega)^d.
 \end{align*}
 And for all $\textbf{v} \in H^{-1}(\Omega)^d$, $(S_h(\textbf{v}),r_h)\in \cX_h\times\cM_h$ is the weak solution of the following discrete problem:
 \begin{align*}
a_f(\widehat{\nu};S_h(\textbf{v}),\textbf{v}_h)+d(\textbf{v},r_h)&=(\textbf{v},\textbf{v}_h),\ &&\forall\, \textbf{v}_h\in\cX_h,\\
d(S_h(\textbf{v}),w_h)&=0,\ &&\forall\, w_h\in\cM_h.
 \end{align*}
 \begin{lemma} For $\textbf{v} \in \cX_{0h}$, there exists constant $C>0$ independent of $h$ such that
\begin{align}\label{eq71}
 \|\textbf{v}\|_{*,h}\leq C\|\nabla S_h(\textbf{v})\|.
 \end{align}
\end{lemma}
Proof: The proof follows the similar lines as in \cite{Caishen2018}. Here we skip it.
 $$\eqno\Box$$
\begin{theorem} Suppose that assumptions  of \textbf{Theorem 5} hold.
We have the following estimates
\begin{align*}
\|\nabla S_h(e^n_\textbf{u}-e^{n-1}_\textbf{u})\|&\leq C\sqrt{\Delta t}(\Delta t+h^{k+1}), \\
\Bigl(\Delta t\sum_{n=1}^{N}\|\nabla S_h(e^n_\textbf{u}-e^{n-1}_\textbf{u})\|^2\Bigr)^{1/2} &\leq C{\Delta t}(\Delta t+h^{k+1}).
 \end{align*}
\end{theorem}
\noindent \textit{Proof:}\quad Using \textbf{Assumption A1} and \eqref{eq37.1a}, we obtain
\begin{align}\label{eq72}
(\delta_t\epsilon_\textbf{u}^n,\textbf{v}_h)+a_f(\nu^n(\phi^{n-1}_h);\epsilon_\textbf{u}^n,\textbf{v}_h)
-d(\textbf{v}_h,\epsilon_p^n)&=(\Phi_h^n,\textbf{v}_h)+(R^n_\textbf{u},\textbf{v}_h),\\\label{eq73}
d(\epsilon_\textbf{u}^n,q_h)&=0.
 \end{align}
Setting $\textbf{v}_h=S_h(\epsilon_\textbf{u}^n-\epsilon_\textbf{u}^{n-1})\in \cX_h$ in \eqref{eq72}, noting that $d\bigl(S_h(\epsilon_\textbf{u}^n)-S_h(\epsilon_\textbf{u}^{n-1}),q_h\bigr)=0$, one finds that
\begin{align}\label{eq74}
\frac{\|\nabla S_h(\epsilon^n_\textbf{u}-\epsilon^{n-1}_\textbf{u})\|^2}{\Delta t}&=-a_f\bigl(\nu^n(\theta^{n-1}_h);\epsilon_\textbf{u}^n,S_h(\epsilon_\textbf{u}^n-\epsilon_\textbf{u}^{n-1})\bigr)\\&\nonumber\quad
+\bigl(R^n_\textbf{u},S_h(\epsilon_\textbf{u}^n-\epsilon_\textbf{u}^{n-1})\bigr)+(\Phi_h^n,S_h(\epsilon_\textbf{u}^n-\epsilon_\textbf{u}^{n-1}).
 \end{align}
Applying the same arguments as those in the proof of \textbf{Theorem 5}, the RHS of \eqref{eq74} can be estimated as
 \begin{align*}
-a_f\bigl(\nu^n(\theta^{n-1}_h);\epsilon_\textbf{u}^n,S_h(\epsilon_\textbf{u}^n-\epsilon_\textbf{u}^{n-1})\bigr)&\leq  \frac{1}{6\Delta t}\|\nabla S_h(\epsilon_\textbf{u}^n-\epsilon_\textbf{u}^{n-1})\|^2+C\Delta t\|\nabla\epsilon_\textbf{u}^n\|^2,\\
\bigl(R^n_\textbf{u},S_h(\epsilon_\textbf{u}^n-\epsilon_\textbf{u}^{n-1})\bigr)&\leq \frac{1}{6\Delta t}\|\nabla S_h(\epsilon_\textbf{u}^n-\epsilon_\textbf{u}^{n-1})\|^2
+C\Delta t(\Delta t^2+h^{2k+2}),\\
(\Phi_h^n,S_h(\epsilon_\textbf{u}^n-\epsilon_\textbf{u}^{n-1})&\leq  \frac{1}{6\Delta t}\|\nabla S_h(\epsilon_\textbf{u}^n-\epsilon_\textbf{u}^{n-1})\|^2
+C\Delta t(\Delta t^2+h^{2k+2})\\&\quad+C\Delta t\bigl(\|\nabla\epsilon_\mu^n\|^2+\|\nabla\epsilon_\textbf{u}^n\|^2
 +\|\nabla\epsilon_\phi^{n-1}\|^2\\&\quad+\|\nabla\epsilon_\textbf{u}^{n-1}\|^2+\|\nabla\epsilon_\textbf{B}^{n}\|^2+\|\nabla\epsilon_\textbf{B}^{n-1}\|^2
\bigr).
 \end{align*}
Combining \eqref{eq74} with above bounds, it follows that
 \begin{align*}
 \frac{\|\nabla S_h(\epsilon^n_\textbf{u}-\epsilon^{n-1}_\textbf{u})\|^2}{2\Delta t}
 &\leq C\Bigr\{\Delta t(\Delta t^2+h^{2k+2})+C\Delta t\bigl(\|\nabla\epsilon_\mu^n\|^2+\|\nabla\epsilon_\textbf{u}^n\|^2\\&\quad
 +\|\nabla\epsilon_\phi^{n-1}\|^2+\|\nabla\epsilon_\textbf{u}^{n-1}\|^2+\|\nabla\epsilon_\textbf{B}^{n}\|^2+\|\nabla\epsilon_\textbf{B}^{n-1}\|^2
\bigr)\Bigr\}.
 \end{align*}
 Applying \textbf{Theorems 5} and the triangle inequality, the desired result is obtained. The proof is finished.
 $$\eqno\Box$$

Then we give and prove the second main result of this
section for the pressure.

\begin{theorem} Suppose that assumptions of \textbf{Theorem 5} and $\Delta t\leq Ch$ hold.
Then the finite element approximate $p^n_h$ in \eqref{eq17.1} satisfy following bound:
\begin{align*}
\Bigl(\Delta t\sum_{n=1}^{N}\|p(t_n)-p^n_h\|^2\Bigr)^{1/2} &\leq C\bigl(\Delta t+h^{k+1}\bigr).
 \end{align*}
\end{theorem}
\noindent \textit{Proof:}\quad From \textbf{Assumption A1} and \eqref{eq37.1a}, one finds that
\begin{align}\label{eq75}
{\beta}_0\|\epsilon^n_p\|&\leq \sup_{\textbf{v}_h\in \cX_h, \textbf{v}_h\neq\textbf{0}}\dfrac{d(\epsilon^n_p,\textbf{v}_h)}{\|\nabla\textbf{v}_h\|}\\\nonumber\quad&
\leq \sup_{\textbf{v}_h\in \cX_h, \textbf{v}_h\neq\textbf{0}}\Bigl\{\dfrac{(\delta_t\epsilon_\textbf{u}^n,\textbf{v}_h)
+a_f(\nu^n(\phi^{n-1}_h);\epsilon_\textbf{u}^n,\textbf{v}_h)}{\|\nabla\textbf{v}_h\|}\Bigr\} \\\nonumber&\quad
+\sup_{\textbf{v}_h\in \cX_h, \textbf{v}_h\neq\textbf{0}}\Bigl\{\dfrac{-(\Phi_h^n,\textbf{v}_h)-(\partial_t\textbf{u}^n-\delta_tR_h{\textbf{u}}^n,\textbf{v}_h)}
{\|\nabla\textbf{v}_h\|}\Bigr\}\\\nonumber
&\leq \sup_{\textbf{v}_h\in \cX_h, \textbf{v}_h\neq\textbf{0}}\dfrac{1}{\|\nabla\textbf{v}_h\|}\Bigl\{b(\textbf{u}_h^{n-1},\epsilon_\textbf{u}^n,\textbf{v}_h)
+S_cc_{\widehat{B}}(\textbf{B}_h^{n-1},\epsilon_\textbf{B}^n,\textbf{v}_h)\\\nonumber&\quad- \lambda(\epsilon_\mu^n\nabla\phi_h^{n-1},\textbf{v}_h)\Bigr\}+C\Bigr\{\|\delta_t\epsilon_\textbf{u}^n\|_{*,h}+\|\nabla\epsilon_\textbf{u}^n\|\\\nonumber&\quad
+ C\bigl(\Delta t+h^{k+1}\bigr)+\bigl(\|\epsilon_{\textbf{u}}^{n-1}\|+\|\epsilon_{\textbf{B}}^{n-1}\|+\|\epsilon_\phi^{n-1}\|\bigr)\Bigr\}\\\nonumber
&\leq \sup_{\textbf{v}_h\in \cX_h, \textbf{v}_h\neq\textbf{0}}\dfrac{1}{\|\nabla\textbf{v}_h\|}\Bigl\{b(\textbf{u}_h^{n-1}-\textbf{u}^{n-1},\epsilon_\textbf{u}^n,\textbf{v}_h)
+b(\textbf{u}^{n-1},\epsilon_\textbf{u}^n,\textbf{v}_h)\\\nonumber&\quad
+S_cc_{\widehat{B}}(\textbf{B}_h^{n-1}-\textbf{B}^{n-1},\epsilon_\textbf{B}^n,\textbf{v}_h)+S_cc_{\widehat{B}}(\textbf{B}^{n-1},\epsilon_\textbf{B}^n,\textbf{v}_h)\\\nonumber&\quad -\lambda\bigl(\epsilon_\mu^n\nabla(\phi_h^{n-1}-\phi^{n-1}),\textbf{v}_h\bigr)-\lambda(\epsilon_\mu^n\nabla\phi^{n-1},\textbf{v}_h)\Bigr\}+C\Bigr\{\|\delta_t\epsilon_\textbf{u}^n\|_{*,h}\\\nonumber&\quad
+\|\nabla\epsilon_\textbf{u}^n\|
+ C\bigl(\Delta t+h^{k+1}\bigr)+\bigl(\|\epsilon_{\textbf{u}}^{n-1}\|+\|\epsilon_{\textbf{B}}^{n-1}\|+\|\epsilon_\phi^{n-1}\|\bigr)\Bigr\}.
 \end{align}
By \textbf{Assumption A2} and \textbf{Theorem 5}, it follows that
\begin{align*}
\|\nabla e_{\textbf{u}}^{n-1}\|&\leq C\min\Bigl\{h^{-1}\|e_{\textbf{u}}^{n-1}\|,\|\nabla e_{\textbf{u}}^{n-1}\|\Bigr\}\\&
\leq C\min\Bigl\{h^{-1}(\Delta t+h^{k+1}), \Delta t^{-1/2}\bigl(\Delta t+h^{k+1}\bigr)\Bigr\}\leq C,\\
\|\nabla e_{\textbf{B}}^{n-1}\|&\leq C\min\Bigl\{h^{-1}\|e_{\textbf{B}}^{n-1}\|,\|\nabla e_{\textbf{B}}^{n-1}\|\Bigr\}\\&
\leq C\min\Bigl\{h^{-1}\bigl(\Delta t+h^{k+1}\bigr), \Delta t^{-1/2}\bigl(\Delta t+h^{k+1}\bigr)\Bigr\}\leq C,\\
\|\nabla e_{\phi}^{n-1}\|&\leq C\min\Bigl\{h^{-1}\|e_{\phi}^{n-1}\|,\|\nabla e_{\phi}^{n-1}\|\Bigr\}\\&
\leq C\min\Bigl\{h^{-1}(\Delta t+h^{k+1}), \Delta t^{-1/2}\bigl(\Delta t+h^{k+1}\bigr)\Bigr\}\leq C.
 \end{align*}
 Using H\"{o}lder inequality and \eqref{eq6.1}, we get
 \begin{align*}
b(\textbf{u}_h^{n-1}-\textbf{u}^{n-1},\epsilon_\textbf{u}^n,\textbf{v}_h)&\leq  C\|\nabla e_\textbf{u}^{n-1}\|\|\nabla\epsilon_{\textbf{u}}^{n}\|\|\nabla\textbf{v}_h\|\leq C\|\nabla\epsilon_\textbf{u}^{n}\|\|\nabla\textbf{v}_h\|,\\
b(\textbf{u}^{n-1},\epsilon_\textbf{u}^n,\textbf{v}_h)&\leq  C\|\nabla\textbf{u}^{n-1}\|\|\nabla\epsilon_{\textbf{u}}^{n}\|\|\nabla\textbf{v}_h\|\leq C\|\nabla\epsilon_\textbf{u}^{n}\|\|\nabla\textbf{v}_h\|,\\
S_cc_{\widehat{B}}(\textbf{B}_h^{n-1}-\textbf{B}^{n-1},\epsilon_\textbf{B}^n,\textbf{v}_h)&\leq  C\|\nabla e_\textbf{B}^{n-1}\|\|\nabla\epsilon_{\textbf{B}}^{n}\|\|\nabla\textbf{v}_h\|\leq C\|\nabla\epsilon_\textbf{B}^{n}\|\|\nabla\textbf{v}_h\|,\\
S_cc_{\widehat{B}}(\textbf{B}^{n-1},\epsilon_\textbf{B}^n,\textbf{v}_h)&\leq C\|\nabla\textbf{B}^{n-1}\|\|\nabla\epsilon_{\textbf{B}}^{n}\|\|\nabla\textbf{v}_h\|
\leq C\|\nabla\epsilon_{\textbf{B}}^{n}\|\|\nabla\textbf{v}_h\|,\\
-\lambda\bigl(\epsilon_\mu^n\nabla(\phi_h^{n-1}-\phi^{n-1}),\textbf{v}_h\bigr)&\leq C\|\epsilon_\mu^n-\overline{\epsilon}_\mu^n\|_{L^4}\|\nabla e_\phi^{n-1}\|\|\nabla\textbf{v}_h\|_{L^4}\leq
C\|\nabla\epsilon_\mu^{n}\|\|\nabla\textbf{v}_h\|.
 \end{align*}
Combining \eqref{eq75} with above inequalities, we derive
 \begin{align*}
 \|\epsilon^n_p\|&\leq C\Bigr\{\|\delta_t\epsilon_\textbf{u}^n\|_{*,h}+\|\nabla\epsilon_\textbf{u}^n\|+\|\nabla\epsilon_\textbf{B}^n\|+\|\nabla\epsilon_\phi^n\|
+ \bigl(\Delta t+h^{k+1}\bigr)\\\nonumber&\quad
+\|\epsilon_{\textbf{u}}^{n-1}\|+\|\epsilon_{\textbf{B}}^{n-1}\|+\|\epsilon_\phi^{n-1}\|\Bigr\}
 \end{align*}
 Due to
 \begin{align*}
\|\delta_t\epsilon_\textbf{u}^n\|_{*,h}&\leq C\Bigl(\frac{1}{\Delta t}\|\nabla S_h(e^n_\textbf{u}-e^{n-1}_\textbf{u})\|\Bigr)\\
&\leq C\frac{\sqrt{\Delta t}}{\Delta t}(\Delta t+h^{k+1}), \\
\Bigl(\Delta t\sum_{n=1}^{N}\|\delta_t\epsilon_\textbf{u}^n\|^2_{*,h}\Bigr)^{1/2} &\leq C\bigl(\Delta t+h^{k+1}\bigr).
 \end{align*}
Using \textbf{Theorems 5, 6} and the triangle inequality, we obtain the desired result. The proof is finished.
 $$\eqno\Box$$

%
%


\section{\label{Sec5} Numerical results }

In this section, some numerical tests are shown to confirm
 the theoretical convergence rates  of the fully discrete scheme \eqref{eq17.1}.
 The experiments have been finished
with applying the finite element libraries from the Fenics Project \cite{Logg2012}.
The following functions $\kappa, \nu, \mu$ are given
 \begin{align*}
\kappa(\phi)=e^{\phi},\ \qquad \nu(\phi)=e^{-\phi},\ \qquad \eta(\phi)=e^{\phi}.
 \end{align*}
 The finite element pair $P_{2}-P_{2}-\textbf{P}_{2}-\textbf{P}_{2}-{P}_{1}$ for the concentration field, the
chemical potential, the velocity field, the magnetic field and the pressure is considered.

  Note that the fully discrete scheme \eqref{eq17.1} is a nonlinear problem, thus we solve
  the problem \eqref{eq17.1} by a Picard type interation.  Namely, we fix the velocity field $\textbf{u}_h$,
  the magnetic field $\textbf{B}_h$ and the pressure $p_h$ at a given time step, then compute for the phase field $\phi_h$ and the
chemical potential $\mu_h$. Then we compute the velocity field $\textbf{u}_h$,
  the magnetic field $\textbf{B}_h$ and the pressure $p_h$ with these updated.

We consider a square domain $\Omega=[0,1]^2$, and the following functions are given
\begin{equation*}
\begin{array}{ll}
\phi(x,y,t)&= 2 + \sin(t)\cos(\pi x)\cos(\pi y),\\
u_1(x,y,t)&= \pi \sin(2\pi y)\sin^2(\pi x)\sin(t),\\
u_2(x,y,t)&=-\pi \sin(2\pi x)\sin^2(\pi y)\sin(t),\\
p(x,y,t)&=\cos(\pi x)\sin(\pi y)\sin(t),\\
B_1(x,y,t)&= \sin(\pi x)\cos(\pi y)\sin(t),\\
B_2(x,y,t)&=-\sin(\pi y)\cos(\pi x)\sin(t),
\end{array}\end{equation*}
as the exact solutions, and some source terms are taken such that
the exact solutions satisfy \eqref{eq2.1}-\eqref{eq4}.

The parameters are set to $S_c= 1$, $\epsilon=0.05$, $\lambda=1$ and $\Delta t=0.1h^2$,
 and the uniform triangles meshes are employed.
We plot the error estimates of the phase field,
the velocity field, the magnetic field
and the pressure
between the numerical solution and the exact solution at $t=0.5$ with different space
sizes in Figure 1. We observe that the rates of
convergence the fully discrete numerical scheme \eqref{eq17.1}
are second order accurate for all variables.
\begin{figure}
\centering
\subfigure[]{
\includegraphics[scale=0.3]{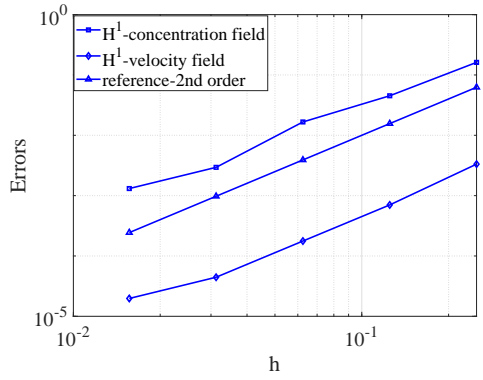}}
\subfigure[]{
\includegraphics[scale=0.3]{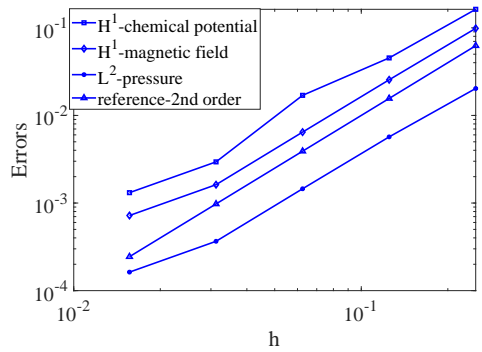}}
\caption{\label{fig:0}\small  Convergence rates: $H^1$ errors of the concentration field $\phi$, the velocity
field $\textbf{u}$, the chemical potential $\mu$, the magnetic field $\textbf{B}$
and $L^2$ the pressure as the space mesh size $h$.}
\end{figure}

\section*{\label{Sec6}Conclusions}

In this paper, we have analyzed a fully discrete  scheme
for computing Cahn-Hilliard-Magneto-hydrodynamics system.
 The scheme is based on using conforming finite element method
in space and Euler semi-implicit discretization with convex splitting in time. We have prove our scheme is unconditionally energy stable and
obtain optimal error estimates for the concentration field, the chemical potential, the velocity field, the magnetic field
and the pressure. Numerical tests are shown to confirm the theoretical rates of the our scheme.


\end{document}